\documentclass[11pt]{amsart}
\usepackage{fullpage}
\usepackage{graphicx}
\usepackage{subfigure}
\usepackage{psfrag}
\usepackage{color}
\usepackage{amsmath}
\usepackage{amssymb}

\newtheorem{theorem}{Theorem}

\newtheorem{remark}[theorem]{Remark}

\title{A numerical study of fluids with pressure dependent viscosity flowing 
  through a rigid porous medium}

\author{K.~B.~Nakshatrala} 
\address{Correspondence to: Dr. Kalyana Babu Nakshatrala, Department of Mechanical Engineering, 
216 Engineering/Physics Building, Texas A\&M University, College Station, Texas - 77843. 
TEL: +1-979-845-1292} \email{knakshatrala@tamu.edu}
\author{K.~R.~Rajagopal}
\address{Professor K.~R.~Rajagopal, Department of Mechanical Engineering, 
Engineering/Physics Building, Texas A\&M University, College Station, Texas - 77843. 
TEL:+1-979-862-4552} \email{krajagopal@tamu.edu}

\date{\today}

\begin{document}

\begin{abstract}
Much of the work on flow through porous media, especially with regard to studies on the 
flow of oil, are based on ``Darcy's law'' or modifications to it such as Darcy-Forchheimer 
or Brinkman models. While many theoretical and numerical studies concerning flow through 
porous media have taken into account the inhomogeneity and anisotropy of the porous solid, 
they have not taken into account the fact that the viscosity of the fluid and drag coefficient 
could depend on the pressure in applications such as enhanced oil recovery. Experiments clearly 
indicate that the viscosity varies exponentially with respect to the pressure and the viscosity 
can change, in some applications, by several orders of magnitude. The fact that the viscosity 
depends on pressure immediately implies that the ``drag coefficient'' would also depend on the 
pressure.   
 
In this paper we consider modifications to Darcy's equation wherein the drag coefficient 
is a function of pressure, which is a realistic model for technological applications like 
enhanced oil recovery and geological carbon sequestration. We first outline the approximations 
behind Darcy's equation and the modifications that we propose to Darcy's equation, and derive 
the governing equations through a systematic approach using mixture theory. We then propose 
a stabilized mixed finite element formulation for the modified Darcy's equation. To solve 
the resulting nonlinear equations we present a solution procedure based on the consistent 
Newton-Raphson method. We solve representative test problems to illustrate the performance 
of the proposed stabilized formulation. One of the objectives of this paper is also to show 
that the dependence of viscosity on the pressure can have a significant effect both on the 
qualitative and quantitative nature of the solution.  
\end{abstract}

\maketitle


\section{INTRODUCTION}
\label{Sec:MDarcy_Introduction}
Flow through porous media is commonly modeled using Darcy's equation \cite{Darcy_1856}. Because 
of its popularity (in many fields like civil, geotechnical, petroleum engineering, composite 
manufacturing) the equation is incorrectly considered as a physical law, and is commonly referred 
to as ``Darcy's law.'' However, it is important to note that Darcy's law is not a new balance law. 
It is an approximation of the balance of linear momentum, and is valid under a plethora of 
assumptions. Some of the serious drawbacks of Darcy's law are
\begin{itemize}
\item it cannot predict the stresses and strains in the solid, 
\item it cannot predict the ``swelling'' and the attendant change in the pore structure, and 
\item it cannot predict the induced inhomogeneity or anisotropy due to the deformation of the solid.
\end{itemize}
In fact, Darcy's law merely predicts the flux of the fluid through the porous solid. More 
importantly, this flux prediction is not accurate at high pressures and pressure gradients, 
which is one of the main themes of this paper. Henceforth, we shall use the terminology 
``Darcy's equation'' instead of ``Darcy's law.'' 

Over the years people have used Darcy's equation beyond its range of applicability. For example, 
Darcy's equation is used to model enhanced oil recovery and geological carbon sequestration. But 
in both these technologically important applications one has to deal with a high pressure range 
of $10 \; \mathrm{MPa}$ - $100 \; \mathrm{MPa}$. There is overwhelming evidence that viscosity 
is not constant and changes drastically with pressure in this range. 

Stokes \cite{Stokes_TCPS_1845_v8_p287} himself recognized that viscosity, in general, depends on 
pressure but that it can be assumed to be constant for a certain class of flows (such as pipe flows 
that involve moderate pressures and pressure gradients). Both enhanced oil recovery and geological 
carbon sequestration do not fall under this class of flows. Barus \cite{Barus_AJS_1893_v45_p87} 
suggested the following exponential relationship between viscosity $\mu$ and pressure $p$ 
\begin{align}
  \label{Eqn:MDarcy_Barus_formula}
  \mu(p) = \mu_0 \exp[\beta p]
\end{align}
where $\beta$ has units $\mathrm{Pa}^{-1}$. Andrade \cite{Bridgman} proposed that the 
variation of viscosity with respect to density ($\rho$), pressure ($p$) and temperature 
($\theta$) can be modeled as
\begin{align}
\label{Eqn:MDarcy_Andrade_formula}
\mu(\rho,p,\theta) = A \rho^{3} \exp\left(\left(p + \rho^2 r\right) \frac{s}{\theta}\right)
\end{align}
where $A$, $r$ and $s$ are empirical constants. A monumental work that discusses in detail 
the effect of pressure on viscosity is the book by Bridgman \cite{Bridgman}, which gives a 
comprehensive account of research concerning the physics of high pressures done prior to 1931. 
It is very important to realize that the flow characteristics of a fluid whose viscosity depends 
on pressure can be significantly different from that of the flow characteristics of a fluid with 
constant viscosity, and this forms one of the main subject matters of this paper. 

We shall now use Barus' formula to get a rough estimate of the variation of the viscosity (and 
hence the drag coefficient) with respect to pressure for some common organic liquids. For Naphthenic 
mineral oil $\beta$ has been determined experimentally to be $23.4 \; \mathrm{GPa}^{-1}$ at $40^{\circ}
\mathrm{C}$ \cite{Hoglund_Wear_1999_v232_p176}. Based on the Barus' formula for Naphthenic mineral oil 
we then have
\begin{align}
  \frac{\mu (p = 100 \; \mathrm{MPa}) - \mu (p = 0.1 \; \mathrm{MPa})}{\mu (p = 0.1 \; \mathrm{MPa})} \times 100 = 935.7 \%
\end{align}
It is worth remarking that since the viscosity varies exponentially, increasing the pressure range 
to that applicable to Elastohydrodynamics would lead to a value of the order of ten to the power 
of eight! On the other hand, using the Dowson-Higginson empirical formula \cite{Dowson_Higginson}, 
the density varies with respect to pressure as 
\begin{align}
  \rho(p) = \rho_0 \left[1 + \frac{0.6 p}{1 + 1.7 p}\right] 
\end{align}
where $p$ is in GPa. Using the Dawson-Higginson formula we have
\begin{align}
  \frac{\rho (p = 100 \; \mathrm{MPa}) - \rho (p = 0.1 \; \mathrm{MPa})}
{\rho (p = 0.1 \; \mathrm{MPa})} = 5.1 \%
\end{align}
Many other experiments have also indicated that change in density due to changes in pressure 
at high pressures is indeed negligible. Thus, it would seem reasonable to neglect the effect 
of pressure on density and consider only the variation of drag function with respect to pressure.

\begin{figure}
  \includegraphics[scale=0.6]{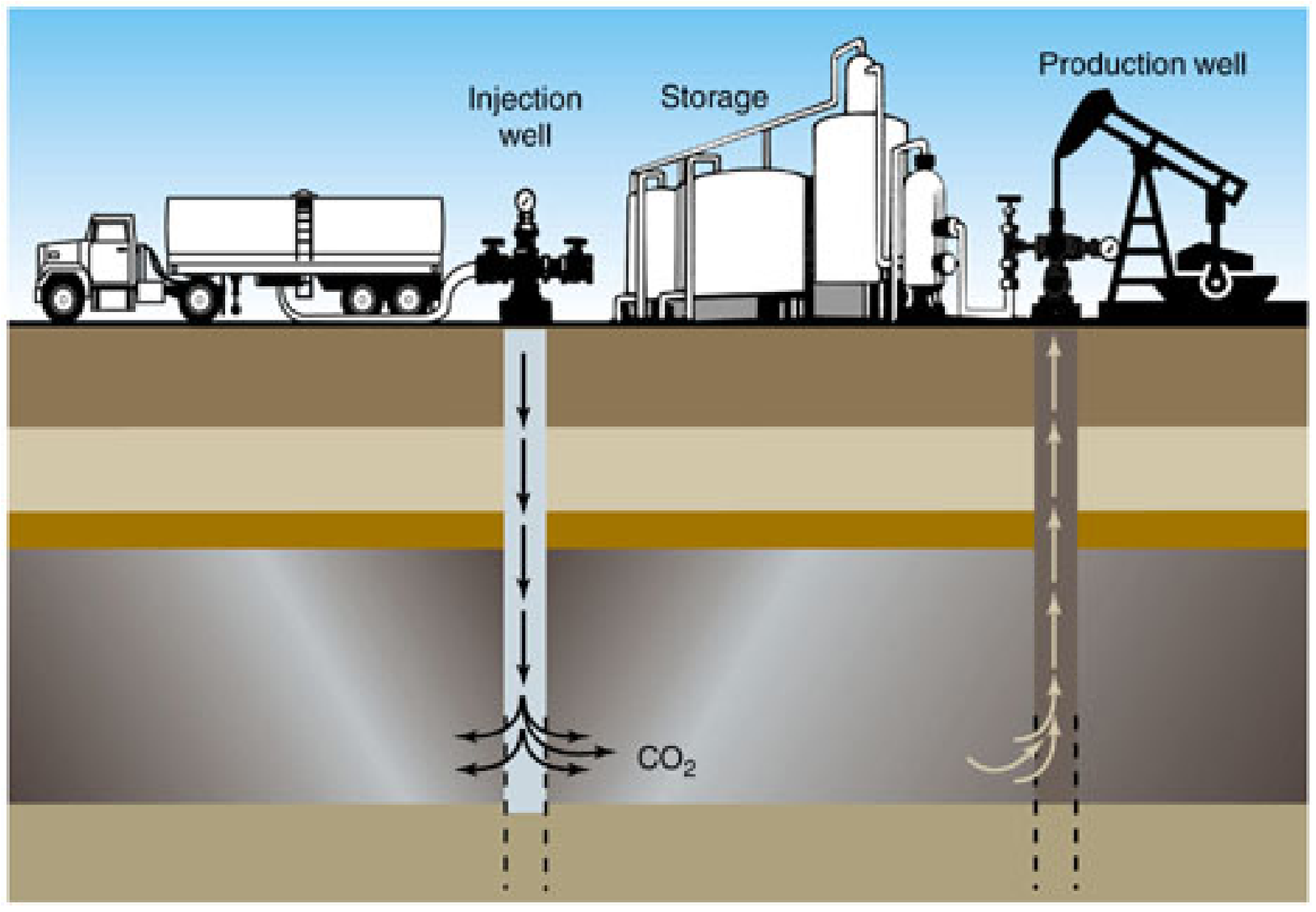}
  \caption{A schematic diagram of enhanced oil recovery. This picture is taken 
    from \textsf{https://www.llnl.gov/str/November01/Kirkendall.html}. 
    \label{Fig:MDarcy_LLNL_EOR}}
\end{figure}

Despite experimental evidence of viscosity depending on pressure in an exponential manner, 
none of the studies of technologically important applications like enhanced oil recovery (EOR) 
and Carbon sequestration take such a dependence on pressure into account. Enhanced oil recovery 
(which is also referred to as tertiary or improved oil recovery) is employed after primary and 
secondary methods of oil recovery, which amount to only 20-40\% of oil recovery. Using EOR techniques 
an additional $30\%$ of oil can be obtained. One of the popular methods of enhanced oil recovery is gas 
injection. In the gas injection method, gas (typically, Carbon dioxide) is pumped into injection 
wells at high pressures, which pushes the oil to the surface at ejection wells. A typical enhanced 
oil recovery process is depicted in Figure \ref{Fig:MDarcy_LLNL_EOR}. In a subsequent section, using 
a typical reservoir simulation, we will show that the results (pressure, pressure gradients, and total 
flow rate) obtained by taking into account the exponential dependence of viscosity on pressure are 
qualitatively and quantitatively different from the classical results. To this end, we consider 
modifying Darcy's equation to take into account the effect of pressure on viscosity. The resulting 
boundary value problem is solved using a stabilized mixed method. 

\subsection{Stabilized mixed finite element formulations}
Numerical methods for Darcy-type equations can be classified into two main categories - primal 
(or single-field) formulations, and mixed (or two-field) formulations. In a primal formulation, 
the whole problem is written in terms of the pressure (which is considered as the primary variable). 
The governing equation will then be a Poisson's equation in terms of the pressure, and one can employ 
any of the standard numerical methods (for example, the Galerkin finite element formulation) to 
solve the resulting equation. Once the pressure field has been calculated, the velocity can be 
calculated using Darcy's equation by using the gradient of the pressure field. 
The two main disadvantages of using a primal formulation have been thoroughly discussed in 
the literature (for example, see \cite{Chen_Huan_Ma,Masud_Hughes_CMAME_2002_v191_p4341}). 
Firstly, for low-order finite elements (which are often employed for generating complex 
meshes) the velocity is poorly approximated as one has to take the gradient of the pressure 
to calculate the velocity. In many situations, velocity is of primary interest. Secondly, 
primal formulations do not possess local mass conservation property with respect to the 
original computational grid. 

On the other hand, mixed formulations can alleviate the drawbacks of primal formulations, 
and are widely used in many numerical simulations. Some of the earlier works on mixed methods 
applied to flows through porous media are \cite{Ewing_Russell_Wheeler_CMAME_1984_v47_p73,
Chavent_Cohen_Jaffre_CMAME_1984_v47_p93,Ewing_Heinemann_CMAME_1984_v47_p161,
Durlofsky_WRR_1994_v30_p965,Arbogast_Wheeler_Yotov_SIAMJNA_1997_v34_p828,
Bergamaschi_Mantica_Manzini_SIAMJSC_1998_v20_p970,Masud_Hughes_CMAME_2002_v191_p4341,
Brezzi_Hughes_Marini_Masud_SIAMJSC_2005_v22_p119,Hughes_Masud_Wan_CMAME_2006_v195_p3347,
Nakshatrala_Turner_Hjelmstad_Masud_CMAME_2006_v195_p4036,Burman_Hansbo_JCAM_2007_v198_p35}. 
However, it is well-known that care should be taken when dealing with mixed formulations. A 
mixed formulation should meet the \emph{inp-sup} (also known as LBB) stability condition 
\cite{Babuska_NumerMath_1971_v16_p322,Brezzi_Fortin}. This falls in the realm of stabilized 
formulations. A huge volume of literature is available on stabilized formulations for various 
equations (e.g., Darcy's equation, Stokes' flow, incompressible Navier-Stokes). A thorough 
discussion of stabilized methods is beyond the scope of this paper. An interested reader 
should refer to \cite{Babuska_Oden_Lee_CMAME_1977_v11_p175,Brooks_Hughes_CMAME_1982_v32_p199,Oden_Jacquotte_CMAME_1984_v43_p231,
Hughes_Franca_Balestra_CMAME_1986_v59_p85,Franca_Hughes_CMAME_1988_v69_p89,Doughlas_Wang_MathComput_1989_v52_p495,
Hughes_Franca_Hulbert_CMAME_1989_v73_p173,Franca_Frey_Hughes_CMAME_1992_v95_p253,
Franca_Russo_CMAME_1997_v142_p361,Masud_Hughes_CMAME_2002_v191_p4341,
Burman_Hansbo_NMPDE_2005_v21_p986,Nakshatrala_Turner_Hjelmstad_Masud_CMAME_2006_v195_p4036,
Correa_Loula_CMAME_2008_v197_p1525} and references therein, and to the following texts concerning 
stabilized methods \cite{Quarternio_Valli_PDE,Roos_Stynes_Tobiska}.

\subsection{Main contributions of this paper} 
Some of the main contributions of this paper are 
\begin{itemize}
\item develop a new stabilized mixed finite element formulation for the modified 
Darcy's equation wherein the drag coefficient (which depends on the viscosity) is 
a function of pressure,
\item document that the dependence of viscosity on the pressure could have a significant 
effect on the nature of the solution (both velocity and pressure fields), and 
\item stress the need and importance of better models than Darcy's equation for 
modeling technologically important problems like enhanced oil recovery and Carbon 
sequestration for which Darcy's equation is physically not appropriate. 
\end{itemize}

\subsection{An outline of the paper} 
The remainder of the paper is organized as follows. In Section \ref{Sec:MDarcy_Assumptions} 
we derive Darcy's equation and our modification of Darcy's equation using mixture theory. In 
Section \ref{Sec:MDarcy_Governing_Equations} we clearly outline our modifications to Darcy's 
equation. In Section \ref{Sec:MDarcy_Stabilized_Formulation} we present a stabilized mixed 
weak formulation for our modification of Darcy's equation, and also present a numerical 
solution procedure based on a consistent Newton-Raphson strategy. Representative numerical 
results are presented in Section \ref{Sec:MDarcy_Numerical_Results}, and conclusions are 
drawn in Section \ref{Sec:MDarcy_Conclusions}. 

%
\section{ASSUMPTIONS BEHIND DARCY'S EQUATION AND MODIFIED DARCY'S EQUATION}
\label{Sec:MDarcy_Assumptions}
We now derive Darcy's equation using mixture theory (which is also referred to as the theory of interacting 
continua) \cite{Atkin_Craine_QJMAM_1976_v29_p209,Bowen,Bedford_Drumheller_IJES_1983_v21_p863,Rajagopal_Tao}. 
This derivation reveals the main assumptions behind Darcy's approximation, and clearly demonstrates the range 
of its (in)applicability. Later, we will also derive a modification of Darcy's equation by relaxing one of the 
assumptions behind Darcy's equation. We first present the balance laws in mixture theory, and then consider the 
flow of a fluid through a (porous) solid. It should be noted that, unless explicitly stated, repeated indices 
do not imply summation.

\subsection{Derivation of modified Darcy equations}
Consider a mixture of $N$ constituents. Let $\rho^{(i)}$ and $\boldsymbol{v}^{(i)}$, respectively, 
denote the density and velocity of the $i$-th constituent. The density of the mixture is defined as 
\begin{align}
\label{Eqn:MDarcy_mixture_density}
\rho = \sum_{i=1}^{N} \rho^{(i)}
\end{align}
The average velocity for the mixture (which is also referred to as the mixture velocity) is defined through 
\begin{align}
\label{Eqn:MDarcy_mixture_velocity}
\boldsymbol{v} = \frac{1}{\rho} \sum_{i=1}^{N} \rho^{(i)} \boldsymbol{v}^{(i)}
\end{align}
It should be noted that a variety of ``mixture velocities" can be defined, and equation 
\eqref{Eqn:MDarcy_mixture_velocity} is one of those possibilities. A detailed discussion 
about this important assumption and its implications with regard to the governing equations 
for mixture can be found in Rajagopal \cite{Rajagopal_M3AS_2007_v17_p215}, and Rajagopal and 
Tao \cite{Rajagopal_Tao}. Let $\boldsymbol{T}^{(i)}$ denote the partial stress of the $i$-th 
constituent. The total stress is defined through 
\begin{align}
\label{Eqn:MDarcy_mixture_total_stress}
\boldsymbol{T} = \sum_{i=1}^{N} \boldsymbol{T}^{(i)}
\end{align}
The total stress can also be defined in several ways, for a discussion of 
the same see Rajagopal and Tao \cite{Rajagopal_Tao}. The balance of mass 
can be written as 
\begin{align}
\label{Eqn:MDarcy_mixture_mass}
\frac{\partial \rho^{(i)}}{\partial t} + \mathrm{div}[\rho^{(i)} \boldsymbol{v}^{(i)}] = m^{(i)} 
\quad \forall i
\end{align}
where $m^{(i)}$ denotes the mass production of the $i$-th constituent. The balance of 
mass for the mixture as a whole warrants that 
\begin{align}
  \label{Eqn:MDarcy_total_mass}
  \sum_{i=1}^{N} m^{(i)} = 0
\end{align}
The balance of linear momentum for the $i$-th constituent can be written as 
\begin{align}
\label{Eqn:MDarcy_mixture_linear_momentum}
\rho^{(i)} \left(\frac{\partial \boldsymbol{v}^{(i)}}{\partial t} + \mathrm{grad}[\boldsymbol{v}^{(i)}] 
\boldsymbol{v}^{(i)}\right) = \mathrm{div}[{\boldsymbol{T}^{(i)}}^{T}] + \rho^{(i)} 
\boldsymbol{b}^{(i)} + \boldsymbol{I}^{(i)} 
\end{align}
where $\boldsymbol{b}^{(i)}$ is the external body force acting on the $i$-th component, and 
$\boldsymbol{I}^{(i)}$ denotes the interaction force that is exerted by other constituents 
on the $i$-th constituent in virtue of their being forced to co-occupy the domain of the 
mixture. By Newton's third law we have 
\begin{align}
\label{Eqn:Mixture_internal_forces}
\sum_{i=1}^{N} \left(\boldsymbol{I}^{(i)} + m^{(i)} \boldsymbol{v}^{(i)}\right) = \boldsymbol{0}
\end{align}
It is important to note that specification of $\boldsymbol{I}^{(i)}$ is also part of 
the constitutive assumptions. In the absence of supply of angular momentum the balance 
of angular momentum can be written as 
\begin{align}
\label{Eqn:Mixture_angular_momentum}
\sum_{i=1}^{N} \boldsymbol{T}^{(i)} = \left(\sum_{i=1}^{N} \boldsymbol{T}^{(i)}\right)^{T}
\end{align}
which basically states that the \emph{total stress} is symmetric. It should be noted that even 
in the absence of angular momentum supply the individual partial stresses $\boldsymbol{T}^{(i)}$ 
need not be symmetric.  

We now consider the case of a mixture consisting of two constituents, namely, a fluid and a solid. 
The quantities associated with the solid and fluid have superscripts `s' and `f', respectively. The 
assumptions behind Darcy's equation and their consequences are as follows:

\begin{enumerate}
\item No mass production of individual constituents. That is, 
  \begin{align}
    m^{(s)} = m^{(f)} = 0 
  \end{align}
\item The solid is assumed to be a rigid porous media. Thus, one can ignore all balance 
  laws for the solid. The stresses in the solid are what they need to be to ensure that 
  the balance of linear momentum in the solid is met.
\item The flow is assumed to be steady 
  \begin{align}
    \frac{\partial \rho^{(f)}}{\partial t} = 0 \quad \mathrm{and} \quad 
    \frac{\partial \boldsymbol{v}^{(f)}}{\partial t} = \boldsymbol{0}
  \end{align}
\item The fluid is assumed to be homogeneous and incompressible. Therefore, 
  \begin{align}
    \mathrm{grad}[\rho^{(f)}] = \boldsymbol{0} \quad \mathrm{and} \quad 
    \mathrm{div}[\boldsymbol{v}^{(f)}] = 0 
  \end{align}
\item The velocity and its gradient are assumed to be small so that the inertial effects can be neglected, 
  which implies that 
  \begin{align}
    \mathrm{grad}[\boldsymbol{v}^{(f)}] \boldsymbol{v}^{(f)} = \boldsymbol{0}
  \end{align} 
\item The partial stress in the fluid is that for an Euler fluid (that is, viscous 
  effects within the fluid are neglected), and the partial stress takes the form 
  \begin{align}
    \boldsymbol{T}^{(f)} = -p^{(f)} \boldsymbol{I}
  \end{align}
\item The only interaction force that comes into play is the frictional force at the 
  boundaries of the pore. It is further assumed that this force is proportional to the 
  relative velocity between the fluid and solid that is reflected by a \emph{drag}-like
  term. Noting the fact that the solid is assumed to be rigid, and attaching the inertial 
  frame to the solid, we have $\boldsymbol{v}^{(s)} = \boldsymbol{0}$. Hence, the 
  interaction force on the fluid takes the form 
  \begin{align}
    \boldsymbol{I}^{(f)} = \alpha \; (\boldsymbol{v}^{(s)} - \boldsymbol{v}^{(f)}) = 
    -\alpha \boldsymbol{v}^{(f)}
  \end{align}
  The drag coefficient $\alpha$, in general, can depend on pressure and relative velocity,  
  $\boldsymbol{v}^{(s)} - \boldsymbol{v}^{(f)}$. 
\item The drag coefficient $\alpha$ is independent of pressure and relative 
  velocity, $\boldsymbol{v}^{(s)} - \boldsymbol{v}^{(f)}$.
\end{enumerate}
\begin{remark}
  In Darcy's equation the drag function is equal to the ratio of the viscosity of 
  the fluid and the coefficient of permeability of the porous medium. That is, 
  \begin{align}
    \alpha = \frac{\mu}{k} 
  \end{align}
  where $k$ is the coefficient of permeability. 
\end{remark}

The above eight assumptions lead to Darcy's equation along with the incompressibility constraint 
\begin{subequations}
\begin{align}
&\alpha \boldsymbol{v}^{(f)} + \mathrm{grad}[p^{(f)}] = \rho \boldsymbol{b}^{(f)} \\
&\mathrm{div}[\boldsymbol{v}^{(f)}] = 0
\end{align}
\end{subequations}

\begin{remark}
One can obtain the Darcy-Forchheimer flow model \cite{Forchheimer_1901_v45_p1782} by relaxing the 
last assumption and allowing the drag coefficient to be a function of the magnitude of the relative 
velocity. Specifically, $\alpha = \alpha_0 + \alpha_1 \|\boldsymbol{v}^{(f)}\|$, where $\alpha_0$ and 
$\alpha_1$ are constants, and $\|\cdot\|$ is the 2-norm (or the Frobenius norm). Note that we have 
used the fact that $\boldsymbol{v}^{(s)} = \boldsymbol{0}$, and hence the magnitude of the relative 
velocity is equal to $\|\boldsymbol{v}^{(f)}\|$. 
\end{remark}

In the next section, we shall relax one of the assumptions behind Darcy's equation by taking into 
account the dependence of viscosity on the pressure, and obtain a modification to Darcy's equation. 
We shall then show that even under this small (but realistic) extension, the response of the fluid is 
dramatically different to the the response of the fluid obtained using Darcy's equation. 
However, it should be notes that the framework offered by theory of interacting continua is quite 
general, and one can accommodate various aspects like mechanical deformation of the solid, 
thermal effects. In our future works, we shall show systematically how the response of various 
components in the mixture change by relaxing various other assumptions under Darcy's equation.

%
\section{GOVERNING EQUATIONS (MODIFIED DARCY EQUATION)}
\label{Sec:MDarcy_Governing_Equations}
The (spatial) position vector is denoted as $\boldsymbol{x}$, and the gradient and divergence operators with 
respect to $\boldsymbol{x}$ are denoted as ``$\mathrm{grad}$" and ``$\mathrm{div}$", respectively. Let $\Omega 
\subset \mathbb{R}^{nd}$ be a bounded open domain, where ``$nd$" denotes the number of spatial dimensions. The 
boundary of $\Omega$ is denoted as $\Gamma$, which is assumed to be piecewise smooth. Mathematically, $\Gamma$ 
is defined as $\Gamma := \bar{\Omega} - \Omega$, where $\bar{\Omega}$ is the closure of $\Omega$. We shall 
denote the velocity vector field as $\boldsymbol{v}(\boldsymbol{x})$. The pressure and density fields are 
denoted as $p(\boldsymbol{x})$ and $\rho(\boldsymbol{x})$, respectively. 

As usual, $\Gamma$ is divided into two parts, denoted by $\Gamma^{\mathrm{D}}$ and $\Gamma^{\mathrm{N}}$, such that 
$\Gamma^{\mathrm{D}} \cap \Gamma^{\mathrm{N}} = \emptyset$ and $\Gamma^{\mathrm{D}} \cup \Gamma^{\mathrm{N}}=\Gamma$. 
$\Gamma^{\mathrm{N}}$ is the part of the boundary on which normal component of the velocity is prescribed, and 
$\Gamma^{\mathrm{D}}$ is part of the boundary on which pressure is prescribed. The modified Darcy equations can 
be written as  
\begin{subequations}
\label{Eqn:MDarcy_modified_Darcy}
\begin{align}
    \label{Eqn:MDarcy_Equilibrium}
    &\alpha(p) \boldsymbol{v} + \mathrm{grad} [p] = \rho \boldsymbol{b} \quad \mbox{in}  \; \Omega \\
    \label{Eqn:MDarcy_Continuity}
    &\mathrm{div}[\boldsymbol{v}] = 0 \quad \mbox{in} \; \Omega \\
    \label{Eqn:MDarcy_Velocity_BC}
    &\boldsymbol{v}(\boldsymbol{x}) \cdot \boldsymbol{n}(\boldsymbol{x}) = v_{n}(\boldsymbol{x}) 
    \quad \mbox{on} \; \Gamma^{\mathrm{N}} \\
    \label{Eqn:MDarcy_Pressure_BC}
    &p(\boldsymbol{x}) = p_0(\boldsymbol{x}) \quad \mbox{on} \; \Gamma^{\mathrm{D}} 
  \end{align}
\end{subequations}
where $\alpha(p)$ (which has dimension of $\left[\mathrm{M} \mathrm{L}^{-3} \mathrm{T}^{-1}\right]$) 
is the drag function, $p_0(\boldsymbol{x})$ is the prescribed pressure, $v_{n}(\boldsymbol{x})$ is 
the prescribed normal component of the velocity, $\boldsymbol{b}(\boldsymbol{x})$ is the specific 
body force, and $\boldsymbol{n}(\boldsymbol{x})$ is the unit outward normal vector to $\Gamma$. 

\begin{remark}
  If $\Gamma^{\mathrm{N}} = \Gamma$ then for well-posedness of \eqref{Eqn:MDarcy_modified_Darcy} 
  one needs to satisfy the following compatibility condition 
  \begin{align}
    \int_{\Gamma^{\mathrm{N}} = \Gamma} v_{n} \; \mathrm{d} \Gamma = 0 
  \end{align}
  which is a direct consequence of the divergence theorem. Also in this case, the solution for 
  pressure is not unique. For uniqueness of the solution one needs to impose an additional 
  constraint on the pressure. For example,  
  \begin{align}
    \label{Eqn:MDarcy_integral_p}
    \int_{\Omega} p \; \mathrm{d} \Omega = 0
  \end{align}
  In numerical simulations, it is common to prescribe the pressure at a point, which is computationally 
  easier than enforcing the constraint \eqref{Eqn:MDarcy_integral_p}. 
\end{remark}

In this paper we consider the following three different drag functions 
\begin{subequations}
  \label{Eqn:MDarcy_specific_alpha}
  \begin{align}
    \label{Eqn:MDarcy_specific_alpha1}
    &\alpha(p) = \alpha_0 \\
    \label{Eqn:MDarcy_specific_alpha2}
    &\alpha(p) = \alpha_0 \left(1 + \beta p \right) \\
    \label{Eqn:MDarcy_specific_alpha3}
    &\alpha(p) = \alpha_0 \exp\left[\beta p\right] 
  \end{align}
\end{subequations}
where $\alpha_0 > 0$ and $\beta > 0$ are constants independent of $p$, $\boldsymbol{v}$ and $\rho$ but 
can be functions of $\boldsymbol{x}$. The dimension of the coefficient $\beta$ is $\left[\mathrm{M}^{-1} 
  \mathrm{L} \mathrm{T}^{-2}\right]$. The constant drag function \eqref{Eqn:MDarcy_specific_alpha1} gives 
rise to (standard) Darcy's equation, and the drag function \eqref{Eqn:MDarcy_specific_alpha2} is based on 
Barus' formula \eqref{Eqn:MDarcy_Barus_formula}. Note that as $p \rightarrow \infty$ the drag coefficient 
\eqref{Eqn:MDarcy_specific_alpha1} $\frac{\alpha(p)}{p} \rightarrow 0$, drag coefficient 
\eqref{Eqn:MDarcy_specific_alpha2} $\frac{\alpha(p)}{p} \rightarrow \alpha_0 \beta$ (a constant), 
and drag coefficient \eqref{Eqn:MDarcy_specific_alpha3} $\frac{\alpha(p)}{p} \rightarrow \infty$.  

Several rigorous mathematical studies have addressed the equations governing the flows of fluids 
with pressure dependent viscosity. Existence of solutions have been established recently by M\'{a}lek 
et al. \cite{Malek_Necas_Rajagopal_ARMA_2002_v165_p243}, Hron et al. \cite{Hron_Malek_Necas_Rajagopal_MCS_2003_v61_p297}, Franta et al. \cite{Franta_Malek_Rajagopal_PRSAMPES_2005_v461_p651}, 
and Buli\v{c}ek et al. \cite{Bulicek_Malek_Rajagopal_IUMJ_2007_v56_p51}. However, it is important to 
note that all the aforementioned existence results have been established under the condition $\frac{\mu}{p} 
\rightarrow \mathrm{constant}$ as $p \rightarrow \infty$. But, experiments suggest that $\frac{\mu}{p} 
\rightarrow \infty$ as $p \rightarrow \infty$ (for example, Barus' formula). Thus, rigorous results such 
as existence of solutions for fluids with pressure dependent viscosity are open with regard to the type 
of variation of the viscosity with pressure observed in experiments.

It is, in general, not possible to obtain analytical solutions for the above boundary value problem 
\eqref{Eqn:MDarcy_modified_Darcy}. One has to resort to numerical solutions to solve realistic problems 
on complex domains. In this paper we employ the Finite Element Method (FEM), which is a powerful and 
systematic technique for numerically solving partial differential equations. To this end, we present 
the classical mixed formulation, and then present a new stabilized mixed formulation for solving 
modified Darcy equations. 

\subsection{Classical mixed formulation}
Let $L^{2}(\Omega)$ be the space of square integrable (scalar) functions on $\Omega$, and 
$\mathbf{L}^{2}(\Omega)$ be space of square integrable vector fields defined on $\Omega$. 
For convenience, we denote the $L^2$ inner-product over a spatial domain, $K$, as
\begin{subequations}
\label{Eqn:MDarcy_L2_inner_product}
\begin{align}
  (a; b)_K = \int_K a \cdot b \; \mathrm{d} K \quad \forall a, b \in L^{2}(K) \\
  (\boldsymbol{a}; \boldsymbol{b})_K = \int_K \boldsymbol{a} \cdot \boldsymbol{b} \; 
  \mathrm{d} K \quad \forall \boldsymbol{a}, \boldsymbol{b} \in \mathbf{L}^{2}(K)
\end{align}
\end{subequations}
The subscript $K$ will be dropped if $K$ is the whole of $\Omega$, that is $K = \Omega$. The 
(Hilbert) spaces $H^{1}(\Omega)$ and $\mathbf{H}(\mathrm{div},\Omega)$ are, respectively, 
defined as
\begin{subequations}
\begin{align}
  H^{1}(\Omega) &:= \left\{a(\boldsymbol{x}) \in L^{2}(\Omega) \; \big| \; \mathrm{grad}[a] 
    \in \mathbf{L}^{2}(\Omega) \right\} \\
  \mathbf{H}(\mathrm{div},\Omega) &:= \left\{\boldsymbol{a}(\boldsymbol{x}) \in \mathbf{L}^{2}(\Omega) 
    \; \big| \; \mathrm{div}[\boldsymbol{a}] \in L^{2}(\Omega) \right\}
\end{align}
\end{subequations}
Furthermore, we shall define the following function spaces
\begin{subequations}
  \begin{align}
    &\mathcal{V} := \left\{\boldsymbol{v}(\boldsymbol{x}) \in \mathbf{H}(\mathrm{div},\Omega) \; \big| \; 
      \mathrm{trace}\left[\boldsymbol{v}(\boldsymbol{x}) \cdot \boldsymbol{n}(\boldsymbol{x})\right] 
      = v_{n}(\boldsymbol{x}) \; \mathrm{on} \; \Gamma^{\mathrm{N}} \right\} \\
    &\mathcal{W} := \left\{\boldsymbol{w}(\boldsymbol{x}) \in \mathbf{H}(\mathrm{div},\Omega) \; \big| \; 
      \mathrm{trace}\left[\boldsymbol{w}(\boldsymbol{x}) \cdot \boldsymbol{n}(\boldsymbol{x})\right] 
      = 0 \; \mathrm{on} \; \Gamma^{\mathrm{N}} \right\} \\
    &\mathcal{P} \equiv L^{2}(\Omega), \quad \mathcal{Q} \equiv H^{1}(\Omega)
  \end{align}
\end{subequations}
where $\mathrm{trace}[\cdot]$ is the standard trace operator \cite{Brezzi_Fortin}. 

The classical mixed formulation for the modified Darcy equation \eqref{Eqn:MDarcy_modified_Darcy} 
can be written as: Find $\boldsymbol{v}(\boldsymbol{x}) \in \mathcal{V}$ and $p(\boldsymbol{x}) 
\in \mathcal{P}$ such that 
\begin{align}
  \label{Eqn:MDarcy_CM_formulation}
  &(\boldsymbol{w}; \alpha(p) \boldsymbol{v}) - (\mathrm{div}[\boldsymbol{w}];p) + 
  (\boldsymbol{w} \cdot \boldsymbol{n};p_0)_{\Gamma^{\mathrm{D}}} - (q;\mathrm{div}[\boldsymbol{v}]) 
  = (\boldsymbol{w};\rho \boldsymbol{b}) \quad \forall \boldsymbol{w}(\boldsymbol{x}) \in \mathcal{W}, 
  \; q(\boldsymbol{x}) \in \mathcal{P}
\end{align}
Assuming the domain and boundary are sufficiently smooth, for the case of constant drag function 
the above weak formulation \eqref{Eqn:MDarcy_CM_formulation} is well-posed \cite{Brezzi_Fortin}. 
A corresponding finite element formulation can be obtained by choosing suitable approximating 
finite element spaces, which (for a conforming formulation) will be finite dimensional subspaces 
of the underlying function spaces of the weak formulation. Let the finite element function spaces 
for the velocity, the weighting function associated with the velocity, and the pressure be denoted 
by $\mathcal{V}^{h} \subset \mathcal{V}$, $\mathcal{W}^h \subset \mathcal{W}$, and $\mathcal{P}^h 
\subset \mathcal{P}$, respectively. 
A conforming finite element formulation of the classical mixed formulation reads: Find 
$\boldsymbol{v}^h(\boldsymbol{x}) \in \mathcal{V}^h$ and $p^h(\boldsymbol{x}) \in \mathcal{P}^h$ 
such that 
\begin{align}
\label{Eqn:MDarcy_FEM_CM_formulation}
&(\boldsymbol{w}^{h}; \alpha(p^{h}) \boldsymbol{v}^{h}) - (\mathrm{div}[\boldsymbol{w}^{h}];p^{h}) + 
(\boldsymbol{w}^{h} \cdot \boldsymbol{n};p_0)_{\Gamma^{\mathrm{D}}} - 
(q^{h};\mathrm{div}[\boldsymbol{v}^{h}]) = (\boldsymbol{w}^{h};\rho \boldsymbol{b}) \quad 
\forall \boldsymbol{w}^{h}(\boldsymbol{x}) \in \mathcal{W}^{h}, \; q^{h}(\boldsymbol{x}) 
\in \mathcal{P}^{h}
\end{align}

For mixed formulations, the inclusions $\mathcal{V}^h \subset \mathcal{V}$, $\mathcal{W}^h \subset \mathcal{W}$, 
and $\mathcal{P}^h \subset \mathcal{P}$ are themselves not sufficient to produce stable results, and additional 
conditions must be met by the members of these finite element spaces to obtain meaningful numerical results. A 
systematic study of these types of conditions on function spaces to obtain stable numerical results is the main 
theme of \emph{mixed finite elements}. One of the main conditions to be met is the Ladyzhenskaya-Babu\v ska-Brezzi 
(LBB) \emph{inf-sup} stability condition. For further details, see references \cite{Brezzi_Fortin,Gunzburger}.

It is well-known that under the classical mixed formulation many combinations of interpolations 
for velocity and pressure produce spurious oscillations (in the pressure field) when applied to 
problems involving incompressibility as a constraint \cite{Brezzi_Fortin,Gunzburger}. In particular, 
the equal-order interpolation for the velocity and pressure (which is computationally the most 
convenient) does not satisfy the LBB condition, and produces spurious oscillations in the pressure 
field. However, elegant solutions for approximating the velocity and pressure that are stable 
under the classical mixed formulation have been proposed \cite{Raviart_Thomas_MAFEM_1977_p292,
Nedelec_NumerMath_1980_v35_p315,Nedelec_NumerMath_1986_v50_p57,Brezzi_Douglas_Marini_NumerMath_1985_v47_p217,
Brezzi_Douglas_Durran_Marini_NumerMath_1987_v51_p237,Brezzi_Douglas_Fortin_Marini_MMNA_1987_v21_p581}. 
These discrete spaces have been successfully used in many numerical simulations, and good accuracy 
has been attained in both the velocity and pressure fields. But a computer implementation of these 
methods is more involved because of different interpolations for the velocity and pressure.

In the next section we present a new stabilized mixed formulation for modified Darcy equations 
under which the computationally convenient equal-order interpolation for velocity and pressure 
is stable.  

%
\section{A STABILIZED MIXED FORMULATION}
\label{Sec:MDarcy_Stabilized_Formulation}
A variational multiscale mixed formulation has been proposed and studied by Hughes 
and Masud \cite{Masud_Hughes_CMAME_2002_v191_p4341} and Nakshatrala \emph{et al.} 
\cite{Nakshatrala_Turner_Hjelmstad_Masud_CMAME_2006_v195_p4036} for Darcy's equation. 
In reference \cite{Nakshatrala_Turner_Hjelmstad_Masud_CMAME_2006_v195_p4036} it has 
been shown that this variational multiscale formulation for Darcy's equation passes 
three-dimensional patch tests even for distorted meshes, and performs well even on 
complex geometries with unstructured meshes. In this paper we extend this formulation 
to the modifications to Darcy's equation developed in this paper. 

A stabilized mixed formulation for modified Darcy's equation \eqref{Eqn:MDarcy_modified_Darcy}
can be written as: Find $\boldsymbol{v}(\boldsymbol{x}) \in \mathcal{V}$ and $p(\boldsymbol{x}) 
\in \mathcal{Q}$ such that 
\begin{align}
  \label{Eqn:MDarcy_stabilized_formulation}
  \left(\boldsymbol{w};\alpha(p) \boldsymbol{v}\right) &- \left(\mathrm{div}[\boldsymbol{w}];p\right) 
  + \left(\boldsymbol{w}\cdot \boldsymbol{n};p_0\right)_{\Gamma^{\mathrm{D}}} - \left(q;\mathrm{div}
    [\boldsymbol{v}]\right) - \left(\boldsymbol{w};\rho \boldsymbol{b}\right) \nonumber \\ 
  &- \frac{1}{2} \left(\alpha(p) \boldsymbol{w} + \mathrm{grad}[q];
    \alpha^{-1}(p)\left(\alpha(p) \boldsymbol{v} + \mathrm{grad}[p] - \rho \boldsymbol{b}\right)\right) 
  = 0 \quad \forall \boldsymbol{w}(\boldsymbol{x}) \in \mathcal{W}, \; q(\boldsymbol{x}) \in \mathcal{Q}
\end{align}
The terms in the first line of equation \eqref{Eqn:MDarcy_stabilized_formulation} are from 
the classical mixed formulation, and the terms in the second line are the stabilization terms. 
The factor $1/2$ on the second line is the stabilization parameter, which is determined neither 
by mesh-dependent parameters nor on material properties (like viscosity). The stabilization term 
is based on the residual of the modified Darcy's equation \eqref{Eqn:MDarcy_modified_Darcy}. 
Similarly, one may add a residual-based term using the incompressibility constraint, which 
in some cases improves stability (e.g., see \cite{Masud_Hughes_CMAME_2002_v191_p4341}). We 
now solve the above weak formulation \eqref{Eqn:MDarcy_stabilized_formulation} using the 
Finite Element Method. 

\subsection{A finite element approximation}
Let the domain $\Omega$ be decomposed into ``$Nele$'' non-overlapping open element subdomains. That is, 
\begin{align}
  \bar{\Omega} = \bigcup_{e = 1}^{Nele} \bar{\Omega}^{e}
\end{align}
where a superposed bar denotes closure. The boundary of $\Omega^e$ is denoted as $\partial \Omega^{e} 
:= \bar{\Omega}^{e} - \Omega^{e}$. For a non-negative integer $m$, $\mathcal{P}^{m}(\Omega^{e})$ denotes 
the linear vector space spanned by polynomials up to $m$th order defined on the subdomain $\Omega^{e}$. We 
shall define the following finite dimensional vector spaces of $\mathcal{V}$, $\mathcal{W}$ and $\mathcal{Q}$:  
\begin{subequations}
  \begin{align}
    \mathcal{V}^{h} &:= \left\{\boldsymbol{v}^{h}(\boldsymbol{x}) \in \mathcal{V} \; \big| \; \boldsymbol{v}^{h}
      (\boldsymbol{x}) \in \left(C^{0}(\bar{\Omega}) \right)^{nd}, \boldsymbol{v}^{h}(\boldsymbol{x}) \big|_{\Omega^e} 
      \in \left(\mathcal{P}^{k}(\Omega^{e})\right)^{nd}, e = 1, \cdots, Nele \right\} \\
    \mathcal{W}^{h} &:= \left\{\boldsymbol{w}^{h}(\boldsymbol{x}) \in \mathcal{W} \; \big| \; \boldsymbol{w}^{h}
      (\boldsymbol{x}) \in \left(C^{0}(\bar{\Omega}) \right)^{nd}, \boldsymbol{w}^{h}(\boldsymbol{x}) \big|_{\Omega^e} 
      \in \left(\mathcal{P}^{k}(\Omega^{e})\right)^{nd}, e = 1, \cdots, Nele \right\} \\
    \mathcal{Q}^{h} &:= \left\{p^{h}(\boldsymbol{x}) \in \mathcal{Q} \; \big| \; p^{h}(\boldsymbol{x}) 
      \in C^{0}(\bar{\Omega}), p^{h}(\boldsymbol{x}) \big|_{\Omega^e} \in \mathcal{P}^{l}(\Omega^{e}), 
      e = 1, \cdots, Nele \right\} 
  \end{align}
\end{subequations}
where $k$ and $l$ are non-negative integers, and recall that ``$nd$'' denotes the number of spatial dimensions. 
A corresponding finite element formulation can be written as: Find $\boldsymbol{v}^{h}(\boldsymbol{x}) 
\in \mathcal{V}^{h}$ and $p^{h}(\boldsymbol{x}) \in \mathcal{Q}^{h}$ such that
\begin{align}
  \label{Eqn:MDarcy_stabilized_FEM}
  (\boldsymbol{w}^{h};\alpha(p^{h}) \boldsymbol{v}^{h}) - (\mathrm{div}
  [\boldsymbol{w}^{h}];p^{h}) + (\boldsymbol{w}^{h} \cdot \boldsymbol{n};p_0)_{\Gamma^{\mathrm{D}}} 
  - (q^{h};\mathrm{div}[\boldsymbol{v}^{h}]) - (\boldsymbol{w}^{h};\rho \boldsymbol{b}) \nonumber \\
  - \frac{1}{2} \left(\alpha(p^{h}) \boldsymbol{w}^{h} + \mathrm{grad}[q^{h}];
    \alpha^{-1}(p^{h}) (\alpha(p^{h}) \boldsymbol{v}^{h} + \mathrm{grad}[p^{h}] - 
    \rho \boldsymbol{b})\right) = 0 \quad \forall \boldsymbol{w}^{h} \in 
  \mathcal{W}^{h}, \; q^{h} \in \mathcal{Q}^{h}
\end{align}
Due to the presence of the term $\frac{1}{2} (\mathrm{grad}[q^{h}];\alpha^{-1}(p^h) \mathrm{grad}[p^h])$, 
the above formulation circumvents the LBB condition, and hence one can employ equal-order interpolation 
for velocity and pressure, which will be confirmed by numerical results presented in a subsequent section. 

We solve the resulting nonlinear equations based on a consistent Newton-Raphson solution procedure. 
We now derive corresponding residual vector and tangent matrix, which are required for such a procedure. 

\subsection{Residual vector}
We shall define the nodal unknowns for the velocity (denoted as $\hat{\boldsymbol{v}}_e$) 
and pressure (denoted as $\hat{\boldsymbol{p}}_e$) for a given element $\Omega^{e}$ as 
\begin{align}
  \label{Eqn:MDarcy_elemental_vCap_pCap}
  \hat{\boldsymbol{v}}_{e} := 
  \left[\begin{array}{ccc}
      v_{1,1} & \cdots & v_{1,nd} \\
      \vdots  & \ddots & \vdots   \\
      v_{n,1} & \cdots & v_{n,nd}
    \end{array} \right], \quad
  \hat{\boldsymbol{p}}_{e} :=\left[\begin{array}{c}
      p_{1}     \\
      \vdots    \\
      p_{n} 
    \end{array} \right]
\end{align}
where $n$ is the number of nodes in the given element, and $nd$ (as mentioned earlier) 
is the number of spatial dimensions. The solution, $\boldsymbol{v}(\boldsymbol{x})$ and 
$p(\boldsymbol{x})$, over the element $\Omega^{e}$ are interpolated in terms of corresponding 
nodal values (which have a superposed hat) as 
\begin{align}
\label{Eqn:MDarcy_Interpolation_for_solution}
\boldsymbol{v}(\boldsymbol{x}) = \hat{\boldsymbol{v}}^T_{e} \boldsymbol{N}^T(\boldsymbol{x}), \quad 
p(\boldsymbol{x}) = \boldsymbol{N}(\boldsymbol{x}) \hat{\boldsymbol{p}}_{e}
\end{align}
where $\boldsymbol{N}(\boldsymbol{x})$ is a row vector of shape functions. Likewise, the weighting 
functions, $\boldsymbol{w}(\boldsymbol{x})$ and $q(\boldsymbol{x})$, over an element are interpolated 
as 
\begin{align}
\label{Eqn:MDarcy_Interpolation_for_weights} 
\boldsymbol{w}(\boldsymbol{x}) = \hat{\boldsymbol{w}}^T_{e} \boldsymbol{N}^T(\boldsymbol{x}), \quad
q(\boldsymbol{x}) = \boldsymbol{N}(\boldsymbol{x}) \hat{\boldsymbol{q}}_{e}
\end{align}
where the nodal (arbitrary) weights $\hat{\boldsymbol{w}}_{e}$ and $\hat{\boldsymbol{q}}_{e}$ 
are defined similar to $\hat{\boldsymbol{v}}_{e}$ and $\hat{\boldsymbol{p}}_{e}$ (which are 
defined in equation \eqref{Eqn:MDarcy_elemental_vCap_pCap}). One can construct residual vectors 
corresponding to equation \eqref{Eqn:MDarcy_stabilized_formulation} by substituting equations 
\eqref{Eqn:MDarcy_Interpolation_for_solution} into equation \eqref{Eqn:MDarcy_stabilized_formulation}, 
and invoking the arbitrariness of $\hat{\boldsymbol{w}}_e$ and $\hat{\boldsymbol{q}}_{e}$. The element 
residual vectors can be written as 
\begin{align}
\label{Eqn:MDarcy_Residual_v}
\boldsymbol{R}^e_{v}(\hat{\boldsymbol{v}}_{e},\hat{\boldsymbol{p}}_{e}) & := \int_{\Omega^e} (\boldsymbol{N}^T \odot \boldsymbol{I}) 
\alpha(p) \boldsymbol{v} \; \mathrm{d} \Omega -\int_{\Omega^e} \mathrm{vec}\left[\boldsymbol{B}^{T}\right] 
p \; \mathrm{d} \Omega + \int_{\partial \Omega^e \cap \Gamma^{\mathrm{D}}} (\boldsymbol{N}^{T} \odot 
\boldsymbol{I}) \boldsymbol{n} p_0 \; \mathrm{d} \Gamma \nonumber \\
&- \int_{\Omega^e} (\boldsymbol{N}^{T} \odot \boldsymbol{I}) \rho \boldsymbol{b} \; 
\mathrm{d} \Omega -\frac{1}{2} \int_{\Omega^e} (\boldsymbol{N}^T \odot \boldsymbol{I}) 
\left(\alpha(p) \boldsymbol{v} + \mathrm{grad}[p] - \rho \boldsymbol{b}\right) \; \mathrm{d} \Omega \\
\label{Eqn:MDarcy_Residual_p}
\boldsymbol{R}^e_{p}(\hat{\boldsymbol{v}}_e,\hat{\boldsymbol{p}}_e) & := -\int_{\Omega^e} \boldsymbol{N}^{T}  
\mathrm{div}[\boldsymbol{v}] \; \mathrm{d} \Omega - \frac{1}{2} \int_{\Omega^e} 
\boldsymbol{B} \; \alpha^{-1}(p) \left(\alpha(p) \boldsymbol{v} + \mathrm{grad}[p] 
- \rho \boldsymbol{b}\right) \; \mathrm{d} 
\Omega   
\end{align}
where $\boldsymbol{v}$ and $p$ are approximated using equation \eqref{Eqn:MDarcy_Interpolation_for_solution}, 
$\mathrm{vec}[\cdot]$ is an operation that represents a matrix as a vector, $\odot$ is Kronecker product 
\cite{Graham_Kronecker} (see Appendix \ref{Sec:MDarcy_Appendix}), $\boldsymbol{B} := \boldsymbol{DN} 
\boldsymbol{J}^{-1}$, $\boldsymbol{J} := \partial \boldsymbol{x}/\partial \boldsymbol{\xi}$ is the element 
Jacobian matrix, $\boldsymbol{DN}$ represents a matrix of (first) derivatives of element shape functions 
with respect to reference coordinates $\boldsymbol{\xi}$. For example, for a three-node triangular element 
the matrix $\boldsymbol{DN}$ will read 
\begin{align}
\boldsymbol{DN} &:= \left[ \begin{array}{cc} \frac{\partial N_1}{\partial \xi_1} & 
\frac{\partial N_1}{\partial \xi_2} \\\vdots & \vdots \\
\frac{\partial N_3}{\partial \xi_1} & \frac{\partial N_3}{\partial \xi_2} \end{array} \right]
\end{align}

We shall define the global unknown vectors $\hat{\boldsymbol{v}}$ and $\hat{\boldsymbol{p}}$ as 
\begin{align}
  \label{Eqn:MDarcy_global_vCap_pCap}
  \hat{\boldsymbol{v}} := 
  \left[\begin{array}{ccc}
      v_{1,1}  & \cdots & v_{1,nd} \\
      \vdots   & \ddots & \vdots   \\
      v_{nn,1} & \cdots & v_{nn,nd}
    \end{array} \right], \quad
  \hat{\boldsymbol{p}} :=\left[\begin{array}{c}
      p_{1}     \\
      \vdots    \\
      p_{nn} 
    \end{array} \right]
\end{align}
where $nn$ denotes the total number of nodes in the finite element mesh. Global residual vectors 
can be constructed from element residual vectors using the standard assembly procedure 
\begin{align}
\boldsymbol{R}_{v}(\hat{\boldsymbol{v}},\hat{\boldsymbol{p}}) = {\huge \boldsymbol{\mathsf{A}}}_{e = 1}^{Nele} 
\boldsymbol{R}^{e}_{v}(\hat{\boldsymbol{v}}_e,\hat{\boldsymbol{p}}_e), \quad 
\boldsymbol{R}_{p}(\hat{\boldsymbol{v}},\hat{\boldsymbol{p}}) = {\huge \boldsymbol{\mathsf{A}}}_{e = 1}^{Nele} 
\boldsymbol{R}^{e}_{p}(\hat{\boldsymbol{v}}_e,\hat{\boldsymbol{p}}_e)
\end{align}
where ${\huge \boldsymbol{\mathsf{A}}}$ is the standard assembly operator \cite{Zienkiewicz}. 
(Recall that $Nele$ denotes the number of elements.) The finite element solution can be obtained 
by solving simultaneously  
\begin{align}
  \boldsymbol{R}_{v} (\hat{\boldsymbol{v}},\hat{\boldsymbol{p}}) = \boldsymbol{0} \; \mathrm{and} \; 
  \boldsymbol{R}_p (\hat{\boldsymbol{v}},\hat{\boldsymbol{p}}) = \boldsymbol{0}
\end{align}

\subsection{Tangent matrix}
\label{Subsec:MDarcy_Tangent_matrix}
Using a Newton-Raphson type approach, one can obtain the solution in an iterative fashion using 
the following update equation until the residual is under a prescribed tolerance 
\begin{subequations}
\label{Eqn:MDarcy_increments}
\begin{align}
\mathrm{vec}\left[\hat{\boldsymbol{v}}^{(i+1)}\right] = 
\mathrm{vec}\left[\hat{\boldsymbol{v}}^{(i)}\right] + \mathrm{vec}
\left[\Delta \hat{\boldsymbol{v}}^{(i)}\right] 
\end{align}
\begin{align}
\mathrm{vec}\left[\hat{\boldsymbol{p}}^{(i+1)}\right] = 
\mathrm{vec}\left[\hat{\boldsymbol{p}}^{(i)}\right] + 
\mathrm{vec}\left[\Delta \hat{\boldsymbol{p}}^{(i)}\right]
\end{align}
\end{subequations}
The updates at iteration $i$ are calculated from the following system of equations
\begin{align}
\label{Eqn:MDarcy_Newton_step}
\left[ \begin{array}{cc}   
\boldsymbol{K}_{vv} & \boldsymbol{K}_{vp} \\
\boldsymbol{K}_{pv} & \boldsymbol{K}_{pp} \\
\end{array} \right] 
\left\{ \begin{array}{c} 
\mathrm{vec}\left[\Delta \hat{\boldsymbol{v}}^{(i)}\right] \\
\mathrm{vec}\left[\Delta \hat{\boldsymbol{p}}^{(i)}\right] \\
\end{array} \right\} 
= - \left\{ \begin{array}{c} 
\boldsymbol{R}_{v} \\
\boldsymbol{R}_{p} \\ \end{array} \right\}
\end{align}
where 
\begin{align}
\boldsymbol{K}_{vv} = {\huge \boldsymbol{\mathsf{A}}}_{e = 1}^{Nele}\boldsymbol{K}^{e}_{vv}, \; 
\boldsymbol{K}_{vp} = {\huge \boldsymbol{\mathsf{A}}}_{e = 1}^{Nele}\boldsymbol{K}^{e}_{vp}, \;
\boldsymbol{K}_{pv} = {\huge \boldsymbol{\mathsf{A}}}_{e = 1}^{Nele}\boldsymbol{K}^{e}_{pv}, \;
\boldsymbol{K}_{pp} = {\huge \boldsymbol{\mathsf{A}}}_{e = 1}^{Nele}\boldsymbol{K}^{e}_{pp}
\end{align}
and the element matrices are defined as 
\begin{subequations}
\label{Eqn:MDarcy_Tangent_matrices}
\begin{align}
\boldsymbol{K}_{vv}^{e} &= \frac{1}{2} \int_{\Omega^e} \left(\boldsymbol{N}^{T} \odot \boldsymbol{I}
\right) \alpha(p) \left(\boldsymbol{N} \odot \boldsymbol{I} \right) \; \mathrm{d} \Omega \\
\boldsymbol{K}_{vp}^{e} &= \frac{1}{2} \int_{\Omega^e} \left(\boldsymbol{N}^{T} \odot \boldsymbol{I}
\right) \boldsymbol{v} \left(\frac{\mathrm{d} \alpha}{\mathrm{d} p}\right) \boldsymbol{N} \; 
\mathrm{d} \Omega - \int_{\Omega^e} \mathrm{vec}\left[\boldsymbol{B}^{T}\right] \boldsymbol{N} 
\; \mathrm{d}\Omega - \frac{1}{2} \int_{\Omega^e} \left(\boldsymbol{N}^{T} \odot \boldsymbol{I} \right) 
\boldsymbol{B}^{T} \; \mathrm{d} \Omega \\
\boldsymbol{K}_{pv}^{e} &= -\int_{\Omega^{e}} \boldsymbol{N}^{T} \mathrm{vec} 
\left[\boldsymbol{B}^{T}\right]^{T} \; \mathrm{d} \Omega - \frac{1}{2} \int_{\Omega^{e}} 
\boldsymbol{B} \left(\boldsymbol{N} \odot \boldsymbol{I} \right) \; \mathrm{d} \Omega \\
\boldsymbol{K}_{pp}^{e} &= - \frac{1}{2} \int_{\Omega^e} \boldsymbol{B} \alpha^{-1}(p) 
\boldsymbol{B}^{T} \; \mathrm{d} \Omega + \frac{1}{2} \int_{\Omega^e} \boldsymbol{B} \; 
\left(\mathrm{grad}[p] - \rho \boldsymbol{b}\right) \frac{1}{\alpha^2(p)} 
\left(\frac{\mathrm{d} \alpha}{\mathrm{d} p} \right)\boldsymbol{N} \; \mathrm{d} \Omega
\end{align}
\end{subequations}

\begin{remark}
  From equation \eqref{Eqn:MDarcy_Tangent_matrices} it is evident that for the modified Darcy's equation, 
  in general, we have ${\boldsymbol{K}_{pv}^{e}} \neq {\boldsymbol{K}_{vp}^{e}}^{T}$, and the matrix 
  $\boldsymbol{K}_{pp}^{e}$ is not symmetric. But in the case of Darcy's equation (for which $\alpha$ 
  is independent of pressure) we do have $\boldsymbol{K}_{pv}^{e} = {\boldsymbol{K}_{vp}^{e}}^{T}$, 
  and the matrix $\boldsymbol{K}_{pp}^{e}$ is symmetric. It should be noted that symmetry of tangent 
  matrix (or the lack of it) is an important factor to be considered in the selection of solvers for 
  the resulting system of linear equations, and also in the design/selection of preconditioners for 
  iterative solvers, which are commonly employed in large-scale problems. 
\end{remark}
%
\section{REPRESENTATIVE NUMERICAL RESULTS}
\label{Sec:MDarcy_Numerical_Results}
\begin{figure}
  \includegraphics[scale=0.5]{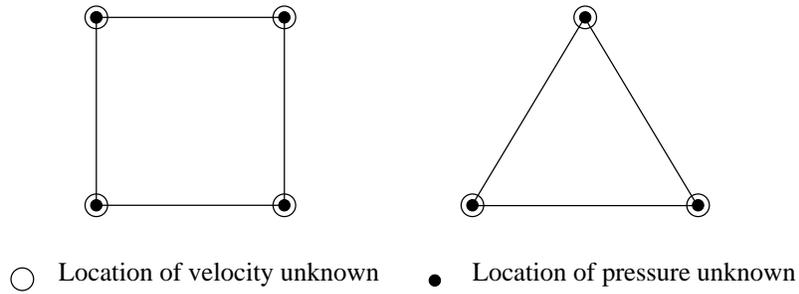}
  \caption{Typical two-dimensional finite elements employed in the numerical simulations, 
    and the location of velocity and pressure unknowns. Equal-order interpolation is used 
    for the velocity and pressure. \label{Fig:MDarcy_elements}} 
\end{figure}

In this section, we present representative numerical results to illustrate the performance of 
the proposed stabilized mixed formulation. In all our numerical simulations we have employed 
equal-order interpolation for pressure and velocity as shown in Figure \ref{Fig:MDarcy_elements}.

First, we shall non-dimensionalize the governing equations. To this end, we define the 
following non-dimensional quantities (which have a superposed bar) 
\begin{align}
  \bar{\boldsymbol{x}} = \frac{\boldsymbol{x}}{L}, \; \bar{\boldsymbol{v}} = \frac{\boldsymbol{v}}{V}, \; 
  \bar{p} = \frac{p}{P}, \; \bar{\alpha} = \frac{\alpha}{\alpha_{\mathrm{ref}}}, \; \bar{\alpha}_0 = 
  \frac{\alpha_0}{\alpha_{\mathrm{ref}}},\; \bar{\rho} = \frac{\rho}{\rho_{\mathrm{ref}}}, \; 
  \bar{\boldsymbol{b}} = \frac{\boldsymbol{b}}{B}, \; \bar{\beta} = \beta P, \; \bar{v}_n = \frac{v_n}{V}, 
  \; \bar{p}_0 = \frac{p_0}{P}
\end{align}
where $L$, $V$, $P$, $\alpha_{\mathrm{ref}}$, $\rho_{\mathrm{ref}}$ and $B$ respectively denote reference length, 
velocity, pressure, drag coefficient, density and specific body force. The gradient and divergence operators with 
respect to $\bar{\boldsymbol{x}}$ are denoted as ``$\overline{\mathrm{grad}}$" and ``$\overline{\mathrm{div}}$", 
respectively. The scaled domain $\Omega_{\mathrm{scaled}}$ is defined as follows: a point in space with position 
vector $\bar{\boldsymbol{x}} \in \Omega_{\mathrm{scaled}}$ corresponds to the same point with position vector given 
by $\boldsymbol{x} = \bar{\boldsymbol{x}} L \in \Omega$. Similarly, one can define the scaled boundaries: $\partial 
\Omega_{\mathrm{scaled}}$, $\Gamma^{\mathrm{D}}_{\mathrm{scaled}}$, and $\Gamma^{\mathrm{N}}_{\mathrm{scaled}}$. 
The above non-dimensionalization gives rise to two dimensionless parameters
\begin{align}
  \mathcal{A} := \frac{\alpha_{\mathrm{ref}} V L}{P}, \quad \mathcal{C}:= \frac{\rho_{\mathrm{ref}} L B}{P}
\end{align}

A corresponding non-dimensionalized form of the drag functions given in equation 
\eqref{Eqn:MDarcy_specific_alpha} can be written as
\begin{subequations} 
  \begin{align}
    \bar{\alpha}(\bar{p}) &= \bar{\alpha}_0 \\ 
    \bar{\alpha}(\bar{p}) &= \bar{\alpha}_0 \left(1 + \bar{\beta} \bar{p}\right) \\
    \bar{\alpha}(\bar{p}) &= \bar{\alpha}_0 \exp[\bar{\beta} \bar{p}] 
  \end{align}
\end{subequations}
We shall write a non-dimensional form of the modified Darcy equation as 
\begin{subequations}
  \begin{align}
    &\mathcal{A} \; \bar{\alpha}(\bar{p}) \bar{\boldsymbol{v}} + \overline{\mathrm{grad}} [p] = 
    \mathcal{C} \; \bar{\rho} \; \bar{\boldsymbol{b}} \quad \; \mathrm{in} \; {\Omega}_{\mathrm{scaled}} \\
    &\overline{\mathrm{div}}[\bar{\boldsymbol{v}}] = 0 \quad \; \mathrm{in} \; {\Omega}_{\mathrm{scaled}} \\
    &\bar{\boldsymbol{v}} \cdot \bar{\boldsymbol{n}} = \bar{v}_n(\bar{\boldsymbol{x}}) \quad \; 
    \mathrm{on} \; {\Gamma}^{\mathrm{N}}_{\mathrm{scaled}} \\
    &\bar{p}(\bar{\boldsymbol{x}}) = \bar{p}_0 (\bar{\boldsymbol{x}}) \quad \; \mathrm{on} \; 
    {\Gamma}^{\mathrm{D}}_{\mathrm{scaled}}
  \end{align}
\end{subequations}
where $\bar{\boldsymbol{n}}(\bar{\boldsymbol{x}})$ is the unit outward normal to the 
boundary $\partial {\Omega}_{\mathrm{scaled}}$. 

\subsection{A simple one-dimensional problem}
\begin{figure}
\psfrag{p1}{$\bar{p}(0) = \bar{p}_1$}
\psfrag{p2}{$\bar{p}(1) = \bar{p}_2$}
\psfrag{x}{$\bar{x}$}
\psfrag{1.0}{$1.0$}
\includegraphics[scale=0.75]{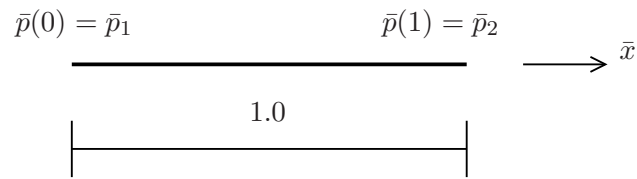}
\caption{Schematic description of the one-dimensional problem. \label{Fig:MDarcy_oneD_problem}}
\end{figure}

We shall test the proposed stabilized mixed formulation using a simple one-dimensional 
problem, which is pictorially described in Figure \ref{Fig:MDarcy_oneD_problem}. Pressures 
of $\bar{p}_1$ and $\bar{p}_2$ are respectively prescribed at the left and right ends of 
the unit domain, and we neglect the body force. The governing equations for this test 
problem can be written as 
\begin{subequations}
  \begin{align}
    &\mathcal{A} \bar{\alpha}(\bar{p}) \bar{v}(\bar{x}) + \frac{\mathrm{d} \bar{p}}{\mathrm{d} \bar{x}} = 0 \quad 
    \mathrm{in} \; (0,1) \\
    &\frac{\mathrm{d} \bar{v}}{d \bar{x}} = 0 \quad \mathrm{in} \; (0,1) \\
    &\bar{p}(\bar{x} = 0) = \bar{p}_1, \quad \bar{p}(\bar{x} = 1) = \bar{p}_2
  \end{align}
\end{subequations}
For this test problem, the analytical solutions with the drag functions defined in equation 
\eqref{Eqn:MDarcy_specific_alpha} can be written as 
\begin{subequations}
  \label{Eqn:MDarcy_OneD_solutions}
  \begin{align}
    \mbox{for the case} \; \bar{\alpha}(\bar{p}) = \bar{\alpha}_0 \; \left\{
      \begin{array}{l}
        \bar{p}(\bar{x}) = \left(\bar{p}_2 - \bar{p}_1 \right) \bar{x} + \bar{p}_1 \\
        \bar{v}(\bar{x}) = - \frac{\left(\bar{p}_2 - \bar{p}_1\right)}{\mathcal{A} \bar{\alpha}_0}
      \end{array} \right.
  \end{align}
  \begin{align}
    \mbox{for the case} \; \bar{\alpha}(\bar{p}) = \bar{\alpha}_0 (1 + \bar{\beta} \bar{p}) \; \left\{
      \begin{array}{l}
        \bar{p}(\bar{x}) = \frac{1}{\bar{\beta}} \left[\left(1 + \bar{\beta}\bar{p}_1\right)^{1 - \bar{x}} 
          \left(1 + \bar{\beta} \bar{p}_2\right)^{\bar{x}} - 1\right] \\
        \bar{v}(\bar{x}) = \frac{-1}{\mathcal{A} \bar{\alpha}_0 \bar{\beta}} \ln 
        \left[\frac{1 + \bar{\beta}\bar{p}_2}{1 + \bar{\beta}\bar{p}_1}\right]
      \end{array} \right.
  \end{align}
  \begin{align}
    \mbox{for the case} \; \bar{\alpha}(\bar{p}) = \bar{\alpha}_0 \exp[\bar{\beta} \bar{p}] \; \left\{
      \begin{array}{l}
        \bar{p}(\bar{x}) = \frac{-1}{\bar{\beta}} \ln \left\{ (1 - \bar{x}) \exp\left[-\bar{\beta}\bar{p}_1\right]  
          + \bar{x} \exp\left[-\bar{\beta} \bar{p}_2\right] \right\} \\
        \bar{v}(\bar{x}) = \frac{1}{\mathcal{A} \bar{\alpha}_0 \bar{\beta}} \left\{
          \exp\left[-\bar{\beta} \bar{p}_2\right] - \exp\left[-\bar{\beta} \bar{p}_1\right] \right\}
      \end{array} \right.
  \end{align}
\end{subequations}

\begin{figure}
  \centering
  \psfrag{p}{$\bar{p}$}
  \psfrag{x}{$\bar{x}$}
  \includegraphics[scale=0.3]{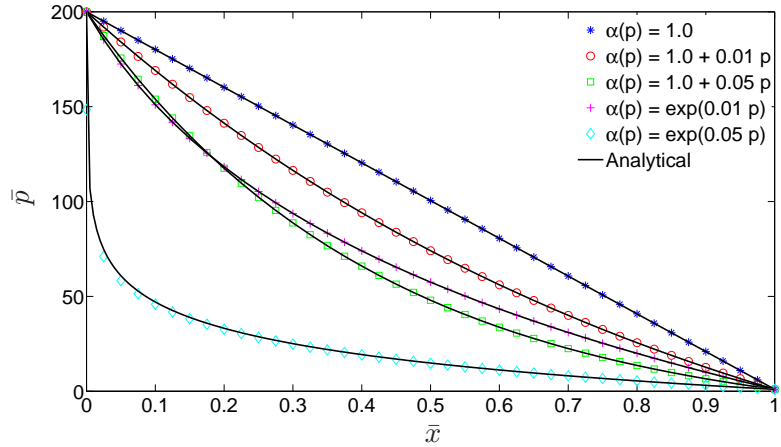}
  \caption{One-dimensional problem: The numerically obtained pressure profiles are compared 
    against analytical solutions for various drag functions. We have chosen $\bar{p}_1 = 200$, 
    $\bar{p}_2 = 1$, $\mathcal{A} = 1$, and $\bar{\alpha}_0 = 1$. The body force is neglected. 
    The test problem is pictorially described in Figure \ref{Fig:MDarcy_oneD_problem}, and the 
    analytical solutions are given in equation \eqref{Eqn:MDarcy_OneD_solutions}. 
    \label{Fig:OneD_problem_results}}
\end{figure}

In Figure \ref{Fig:OneD_problem_results} we compare the numerical solutions against the 
analytical solutions for various drag functions. In all the considered cases, the proposed 
numerical formulation performed well. 

\subsection{Five-spot problem}
\begin{figure}
  \psfrag{vx}{$v_x = 0$}
  \psfrag{vy}{$v_y = 0$}
  \includegraphics[scale=0.75]{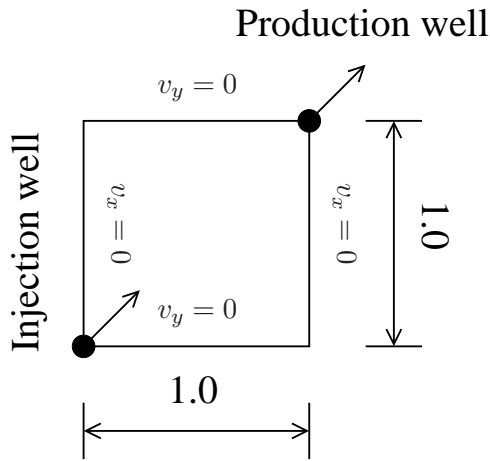}
  \caption{Pictorial description of a quarter of five-spot problem.
    \label{Fig:MDarcy_Five_spot_description}}
\end{figure}

\begin{figure}
  \subfigure[Pressure contours using Darcy's equation]{
    \includegraphics[width=3.5in]{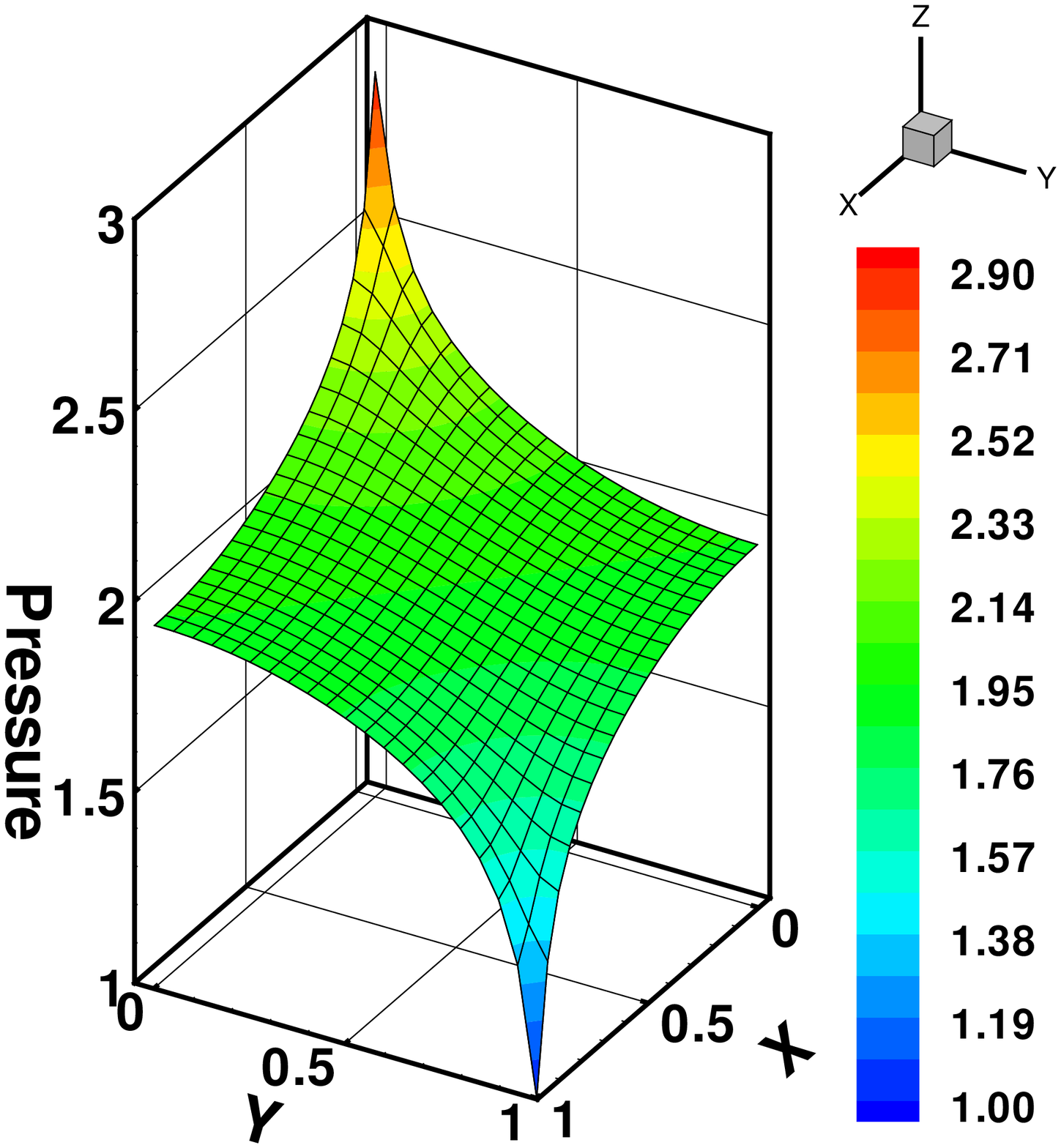}}
  \subfigure[Pressure contours using modified Darcy's equation]{
    \includegraphics[width=3.5in]{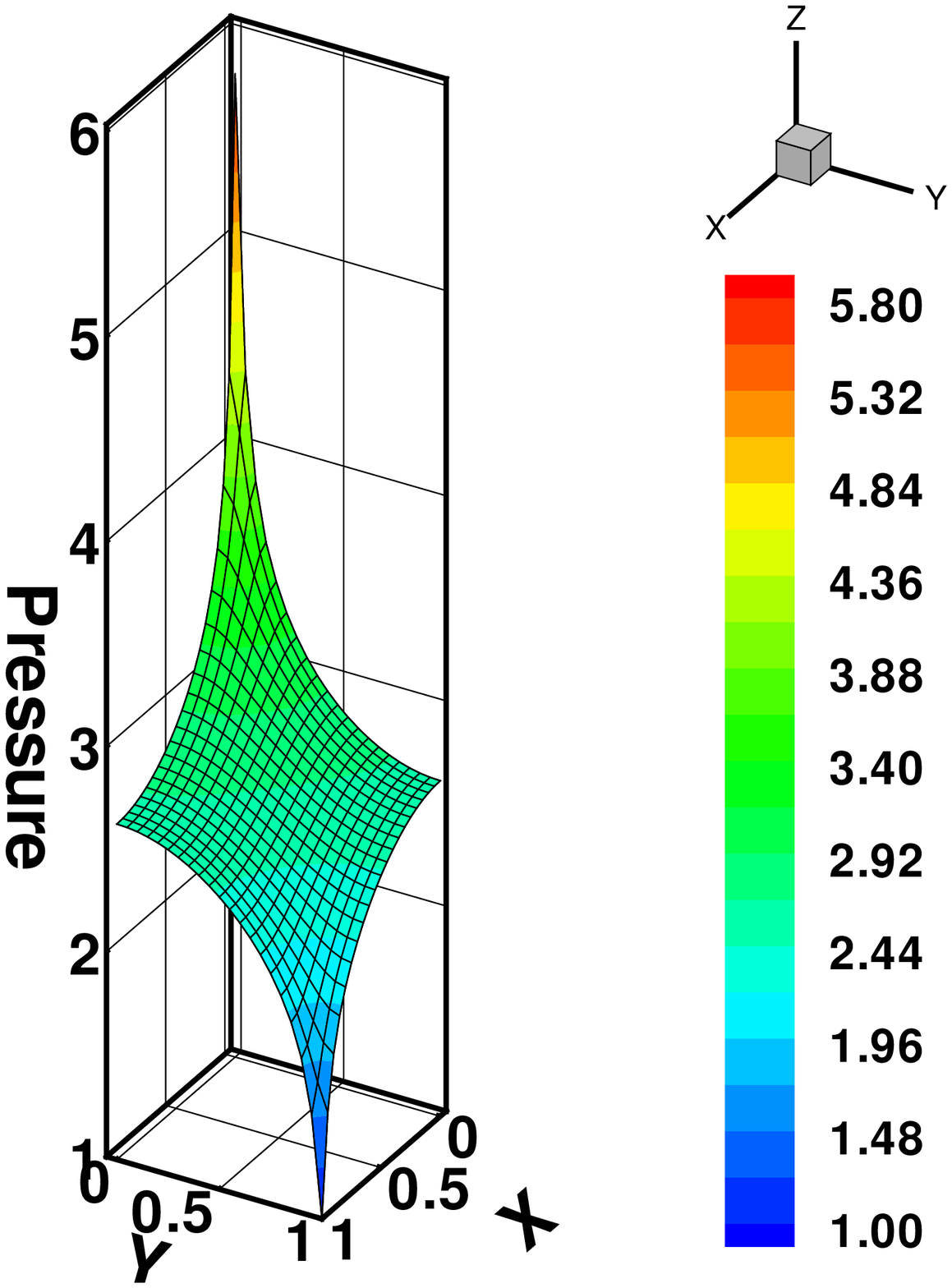}}
  \caption{Five-spot problem using $400$ \emph{four-node quadrilateral elements}: This 
    figure compares elevation plots for the pressure profile using Darcy's equation (top) 
    and modified Darcy's equation (bottom). The pressure based on modified Darcy's equation 
    with the viscosity given by Barus' formula exhibits steeper gradients compared to ones 
    obtained using Darcy's equation. As one can see there are no spurious oscillations in 
    the pressure field despite the steep gradients in pressure near the injection and pressure 
    wells. For Barus' formula we have used $\alpha_0 = 1$ and $\beta = 0.3$. 
    \label{Fig:MDarcy_Five_spot_pressure_profile_Q4}}
\end{figure}

\begin{figure}
  \subfigure[Pressure contours using Darcy's equation]{
    \includegraphics[width=3.5in]{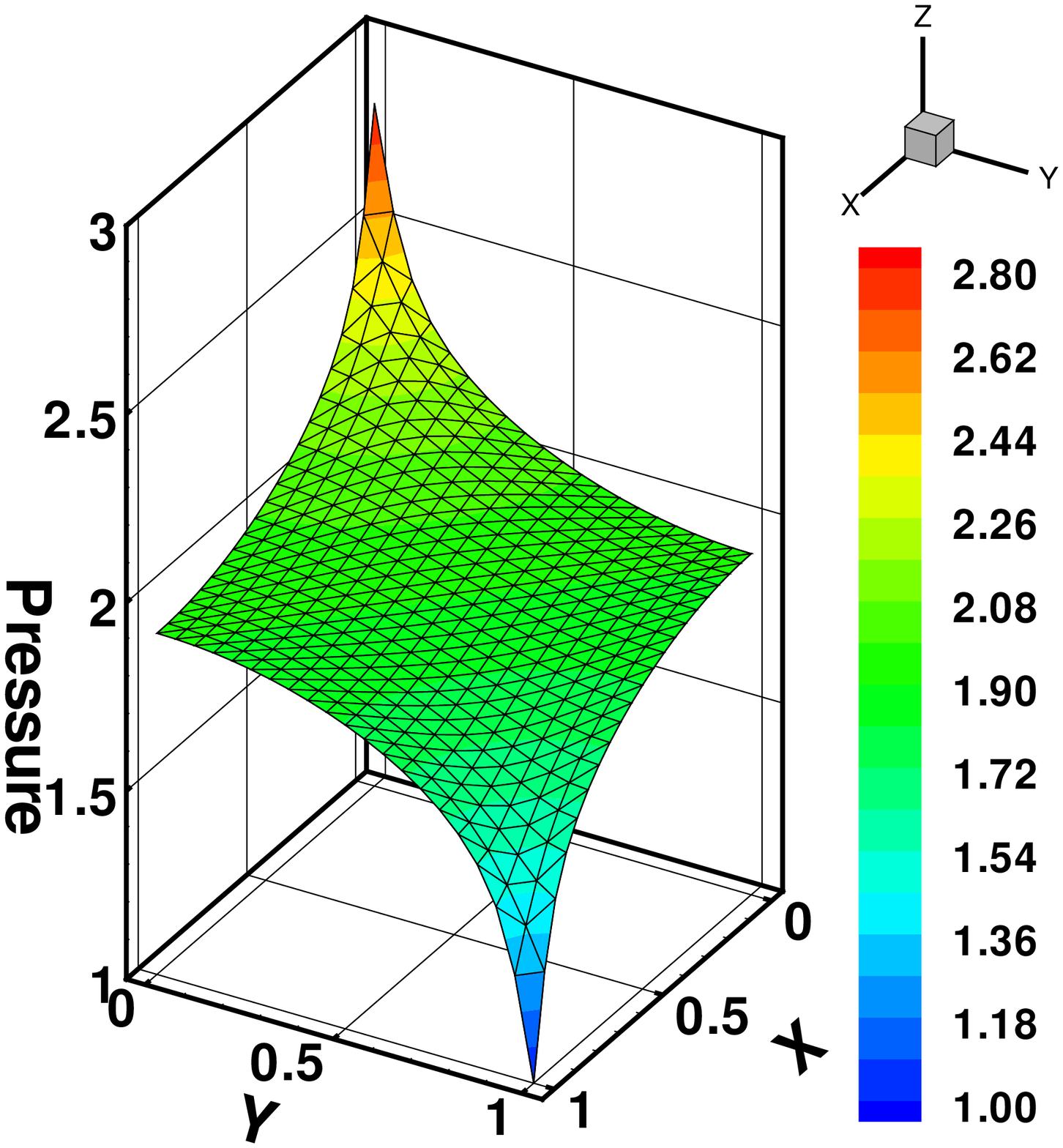}}
  \subfigure[Pressure contous using modified Darcy's equation]{
    \includegraphics[width=3.5in]{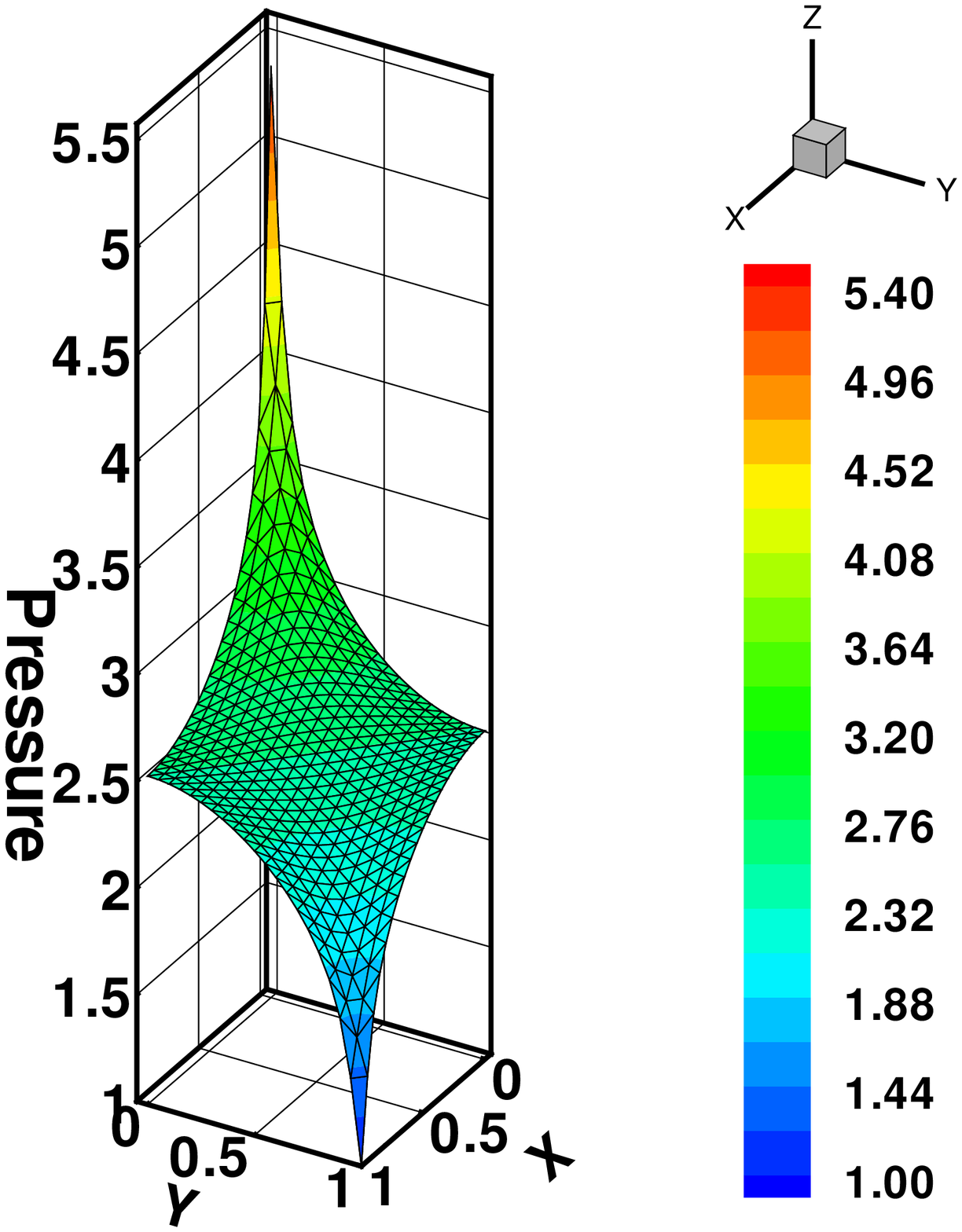}}
  \caption{Five-spot problem using $800$ \emph{three-node triangular elements}: This figure 
    compares elevation plots for the pressure profile using Darcy's equation and modified 
    Darcy's equation (bottom). The pressure based on modified Darcy's equation with the 
    viscosity given by Barus' formula exhibits steeper gradients compared to ones obtained using 
    Darcy's equation. As one can see there are no spurious oscillations in the pressure 
    field despite the steep gradients in pressure near the injection and production wells. 
    For Barus' formula, we have used $\alpha_0 = 1$ and $\beta = 0.3$. 
    \label{Fig:MDarcy_Five_spot_pressure_profile_T3}}
\end{figure}

\begin{figure}
  \subfigure{
    \includegraphics[width=5in]{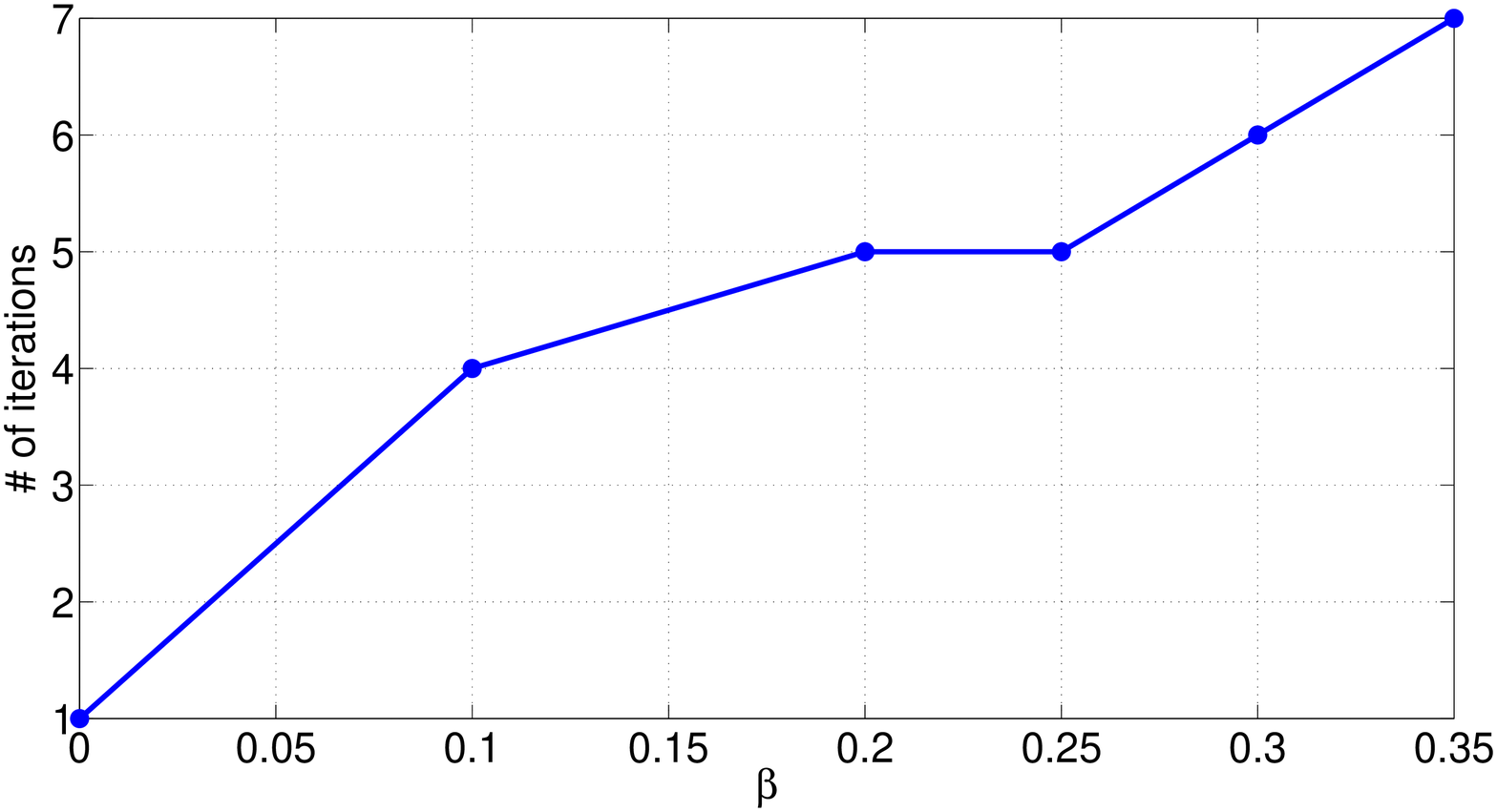}}
  \subfigure{
    \includegraphics[width=5in]{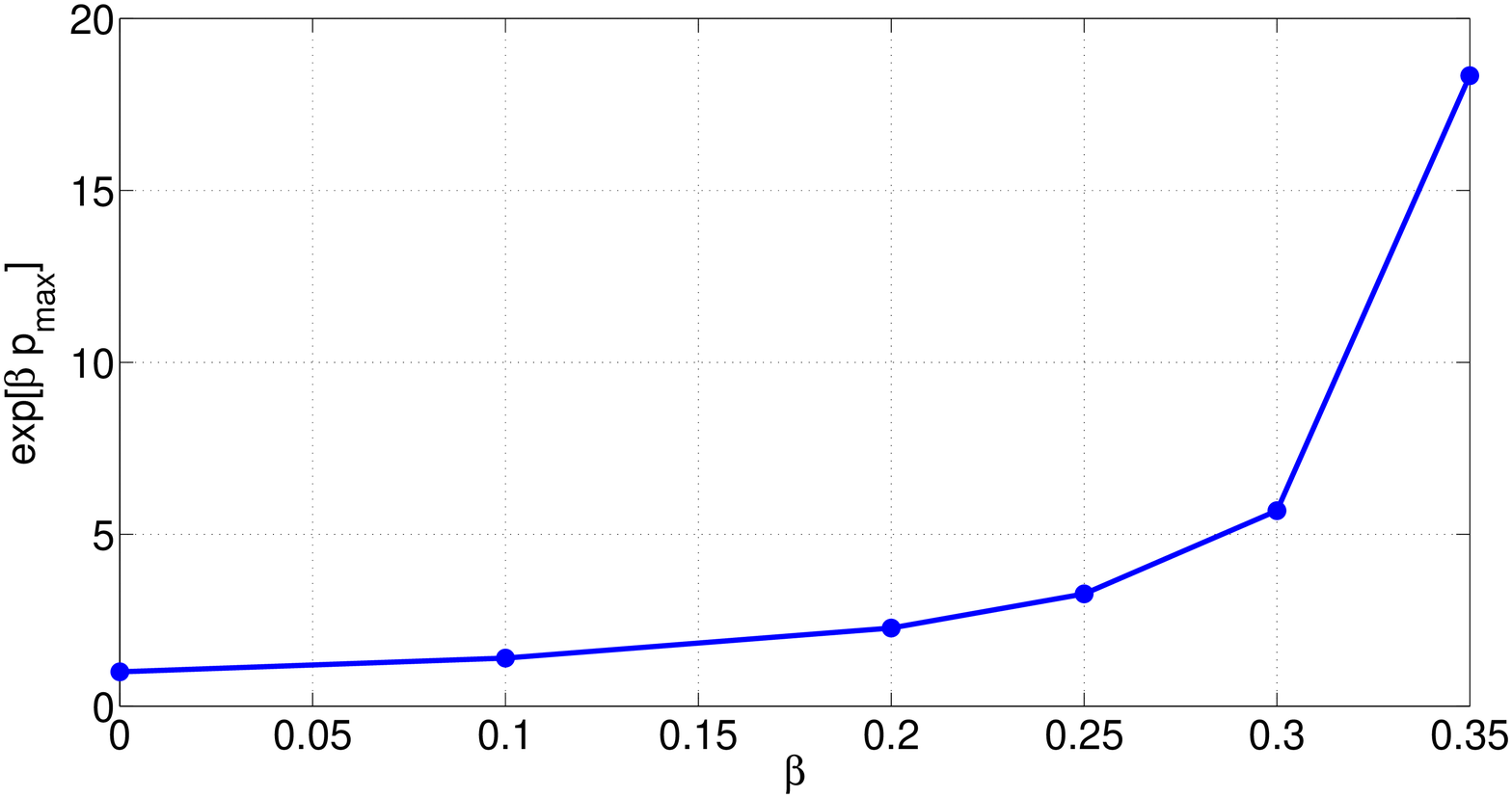}}
  \caption{Five-spot problem using $400$ four-node quadrilateral elements: The top 
  figure shows the number of Newton-Raphson iterations for various values of $\beta$. 
  The bottom figure shows the variation of $\exp[\beta p_{\mathrm{max}}]$ for various 
  values of $\beta$, where $\exp[\beta p_{\mathrm{max}}]$ denotes the maximum pressure 
  in the computational domain.  We have employed Barus' formula with $\alpha_0 = 1$.  
    \label{Fig:MDarcy_Five_spot_different_beta_Q4}}
\end{figure}

This subsection presents numerical results for a quarter of the five-spot problem, which 
exhibits elliptic singularities near (injection and production) wells and is widely used as 
a good model problem to test the robustness of a numerical formulation. Taking into account 
the quarter symmetry in the five-spot problem, the source and sink strengths at injection 
and production wells are, respectively, taken as $+1/4$ and $-1/4$. There is no volumetric 
source/sink (i.e., $\boldsymbol{b}(\boldsymbol{x}) = \boldsymbol{0}$). The five-spot problem 
is pictorially described in Figure \ref{Fig:MDarcy_Five_spot_description}. The numerically 
obtained pressure profiles using Darcy's equation and modified Darcy's equation with the 
viscosity given by Barus' formula are shown in the Figures 
\ref{Fig:MDarcy_Five_spot_pressure_profile_Q4} and \ref{Fig:MDarcy_Five_spot_pressure_profile_T3}. 
Firstly, there are no spurious oscillations in the pressure field even for the modified Darcy's 
equation, and the proposed mixed stabilized formulation performed well. Secondly, for this test 
problem, pressure in the case of the modified Darcy equation exhibits steeper gradients compared 
to the pressure field based on Darcy's equation. 
In Table \ref{Table:MDarcy_NR_convergence} we have shown that the norm of the residual 
vector exhibits terminal quadratic convergence, which basically shows that the tangent 
matrix corresponding to the given residual vector is correctly formulated and calculated. 
In nonlinear problems, (due to the lack of analytical solutions) the terminal quadratic 
convergence is considered as a good measure to check the computer implementation of a 
numerical formulation. 
In Figure \ref{Fig:MDarcy_Five_spot_different_beta_Q4} we compared the number of Newton-Raphson 
iterations for various values of $\beta$ using four-node quadrilateral elements for a relative tolerance of 
$10^{-12}$. In the same figure we have also shown the variation of $\exp[\beta p_{\mathrm{max}}]$ with 
respect to $\beta$, where $p_{\mathrm{max}}$ is the maximum pressure in the computational domain. 
(Note that $p_{\mathrm{max}}$ occurs at the injection well.) As one can see, as $\beta$ increases 
$\exp[\beta p_{\mathrm{max}}]$ increases drastically, and very steep pressure gradients occur near 
the injection and production wells. To obtain results for much higher values of $\beta$ one may have 
to precondition the resulting linear equations at every Newton-Raphson iteration in order to avoid 
ill-conditioned matrices, and also employ a finer mesh.

\begin{remark}
  It is well-known that the Newton-Raphson solution procedure will not always exhibit terminal 
  quadratic convergence. For example, if the multiplicity of the desired root is greater than 
  one, the Newton-Raphson solution procedure may not exhibit terminal quadratic convergence. 
  Also, in some cases, the Newton-Raphson procedure may not even converge. For a detailed 
  discussion of the convergence properties of the Newton-Raphson solution procedure see 
  Reference \cite{Kelley_SIAM_Newtons_method}.
\end{remark}

\begin{table}
  \caption{Five-spot problem: Convergence of Newton-Raphson solution procedure for solving 
    modified Darcy's equation using four-node quadrilateral (Q4) and three-node triangular 
    (T3) elements. We have used Barus' formula with $\alpha_0 = 1$ and $\beta = 0.3$. The 
    test problem is pictorially described in Figure \ref{Fig:MDarcy_Five_spot_description}. 
    As one can see from the table, the solution procedure exhibits terminal quadratic 
    convergence. \label{Table:MDarcy_NR_convergence}}
  \begin{tabular}{ccc} \\ \hline
    Iteration \# &  Norm of residual (Q4) & Norm of residual (T3) \\ \hline           
    0  &   0.242406260    &  0.291915904    \\
    1  &   0.082035285    &  0.086348093    \\
    2  &   0.037068940    &  0.037721275    \\
    3  &   0.004513442    &  0.004318753    \\
    4  &   4.1044706e-05  &  3.4343232e-05  \\
    5  &   1.4244950e-09  &  7.6123308e-10  \\
    6  &   2.0953922e-15  &  3.4935047e-15  \\ \hline
  \end{tabular}
\end{table}

\subsection{$h$-convergence analysis for a two-dimensional problem}
\begin{figure}
  \subfigure[]{
    \includegraphics[scale=0.3]{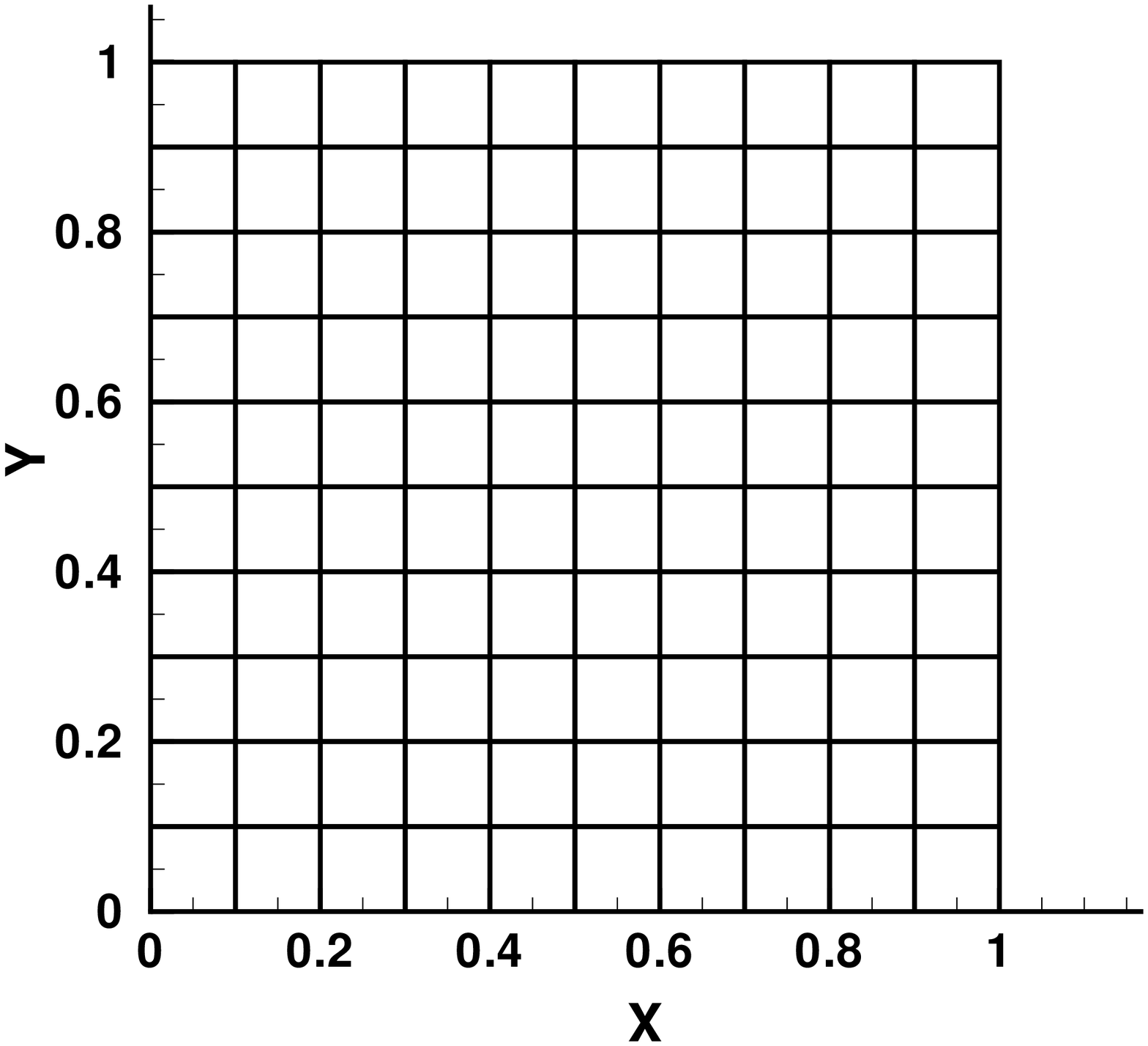}}
  \subfigure[]{
    \includegraphics[scale=0.3]{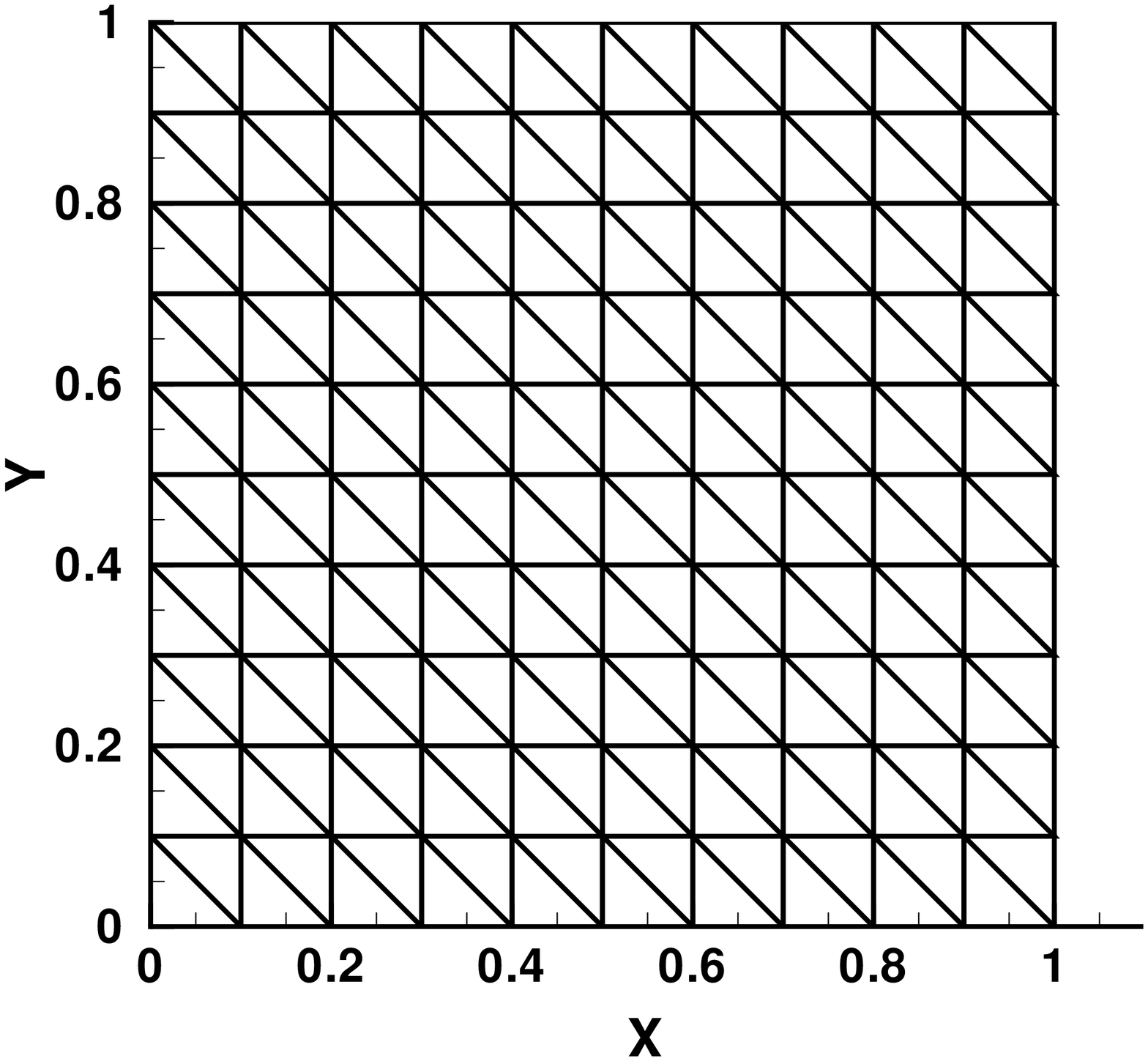}}
  \caption{Typical meshes used in the $h$-convergence analysis: (a) four-node quadrilateral 
    (b) three-node triangular (right) meshes.\label{Fig:MDarcy_FE_convergence_meshes}}
\end{figure}

Let us consider a bi-unit square $[0,1] \times [0,1]$ as our computational domain. 
Let us assume that the velocity and pressure are respectively given as
\begin{align}
  \label{Eqn:2D_problem_exact_solution}
  v_x(x,y) = \sin(\pi x) \cos(\pi y), \quad v_y(x,y) = -\cos(\pi x) \sin(\pi y), \quad p(x,y) = 1 + 25 xy(x - 1)(y - 1) 
\end{align}
We shall assume that $\rho = 1$, $\alpha_0 = 1$, and $\beta = 2$. Using the above assumed solution, 
the required body force under the Barus' formula \eqref{Eqn:MDarcy_Barus_formula} should be equal to 
\begin{subequations}
  \label{Eqn:2D_problem_body_force}
  \begin{align}
    b_x(x,y) &= \exp[\beta (1 + 25 xy (x - 1) (y - 1))] \sin(\pi x) \cos(\pi y) + 25 (2 x - 1) y (y - 1) \\
    b_y(x,y) &= -\exp[\beta (1 + 25 x y (x - 1) (y - 1))] \cos(\pi x) \sin(\pi y) + 25 x (x - 1) (2 y - 1) 
  \end{align}
\end{subequations}
and the boundary conditions are given by 
\begin{align}
  \label{Eqn:2D_problem_BCs}
  v_x(x = 0, y) = v_x(x = 1, y) = 0, \quad v_y(x, y = 0) = v_y(x, y = 1) = 0
\end{align}
Using the aforementioned boundary value problem, we shall perform numerical convergence studies of 
the proposed stabilized mixed formulation for both four-node quadrilateral and three-node triangular 
elements. Typical uniform finite element meshes used in the $h$-convergence analysis are shown in 
Figure \ref{Fig:MDarcy_FE_convergence_meshes}. Note that all the elements in the four-node 
quadrilateral mesh are squares, and all the elements in the three-node triangular mesh are isosceles 
right-angled triangles. In the numerical convergence studies presented in this subsection, we have 
taken $h$ to be equal to the length of the side for square elements, and length of the base (or height) 
for the isosceles right-angled triangles.

\begin{figure}
  \subfigure[Using four-node quadrilateral elements]{
    \includegraphics[scale=0.3]{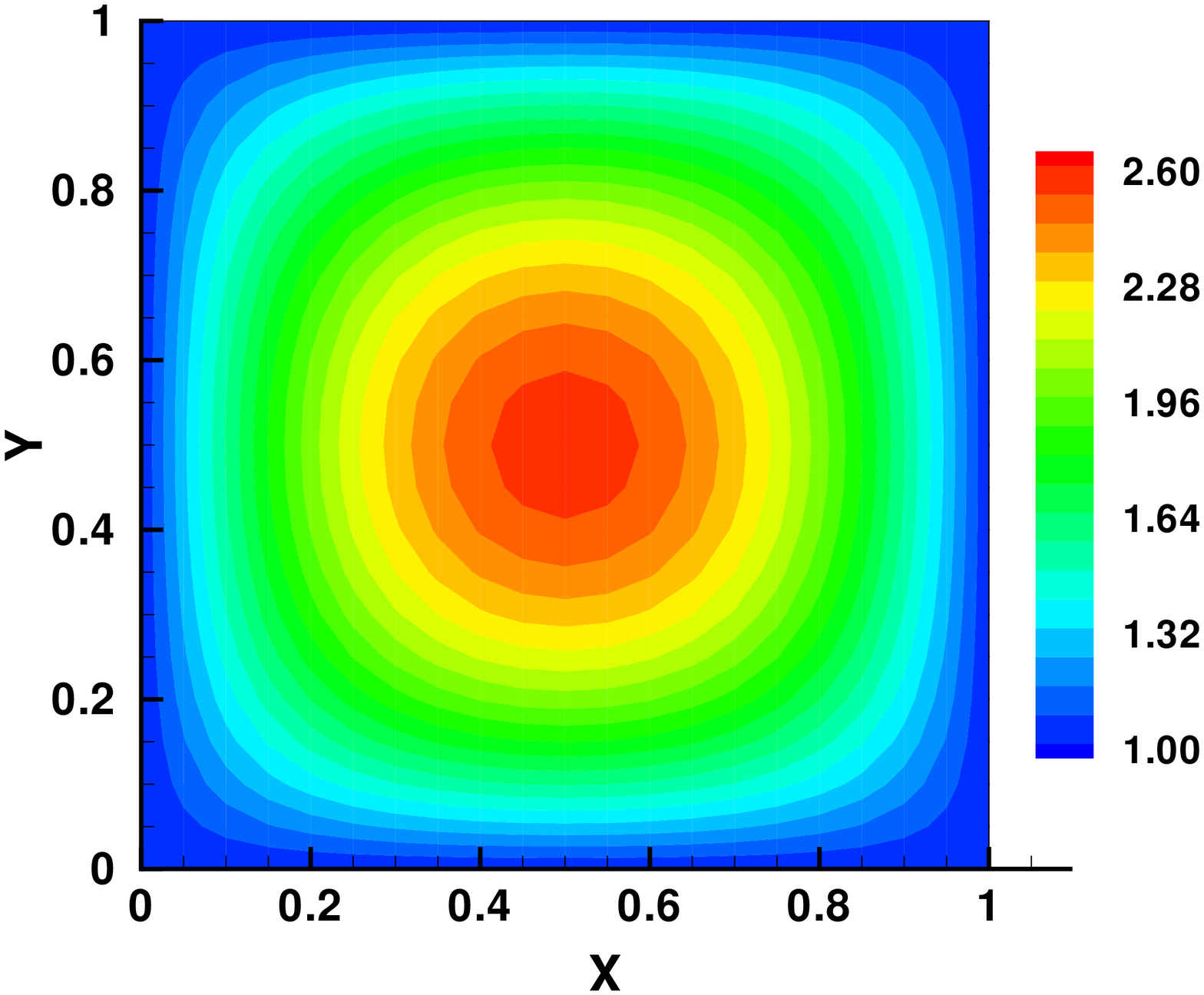}}
  \subfigure[Using three-node triangular elements]{
    \includegraphics[scale=0.3]{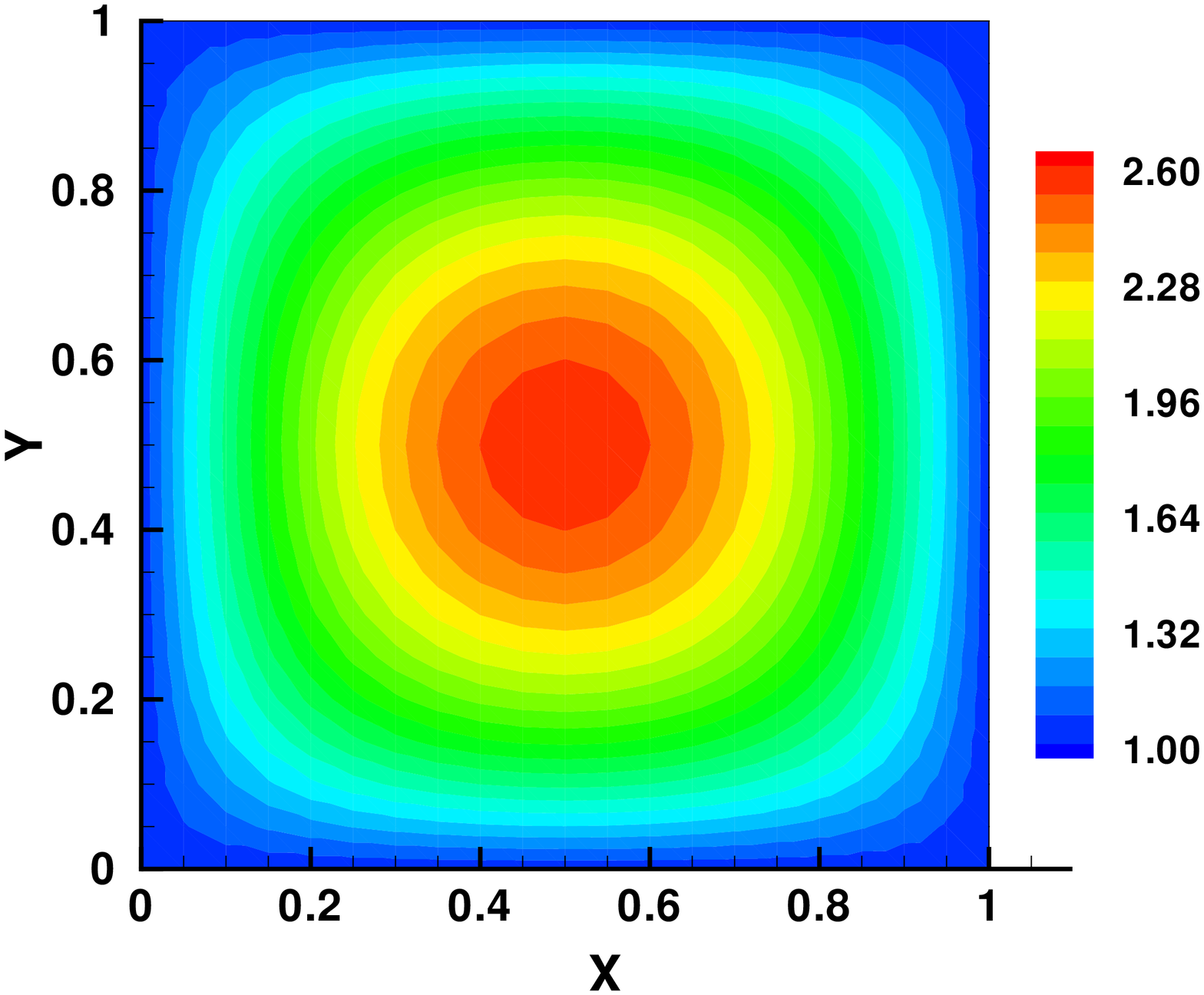}}
  \caption{Structured meshes: contours of pressure using $100$ uniform four-node quadrilateral 
    (left) and $200$ three-node triangular (right) elements. There are no spurious oscillations 
    in the pressure field. \label{Fig:MDarcy_2D_Problem_convergence_Pre}}
\end{figure}

\begin{figure}
  \subfigure[$x$-velocity using Q4 elements]{
    \includegraphics[scale=0.3]{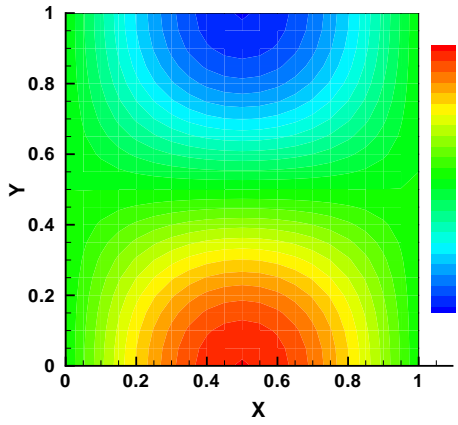}}
  \subfigure[$x$-velocity using T3 elements]{
    \includegraphics[scale=0.3]{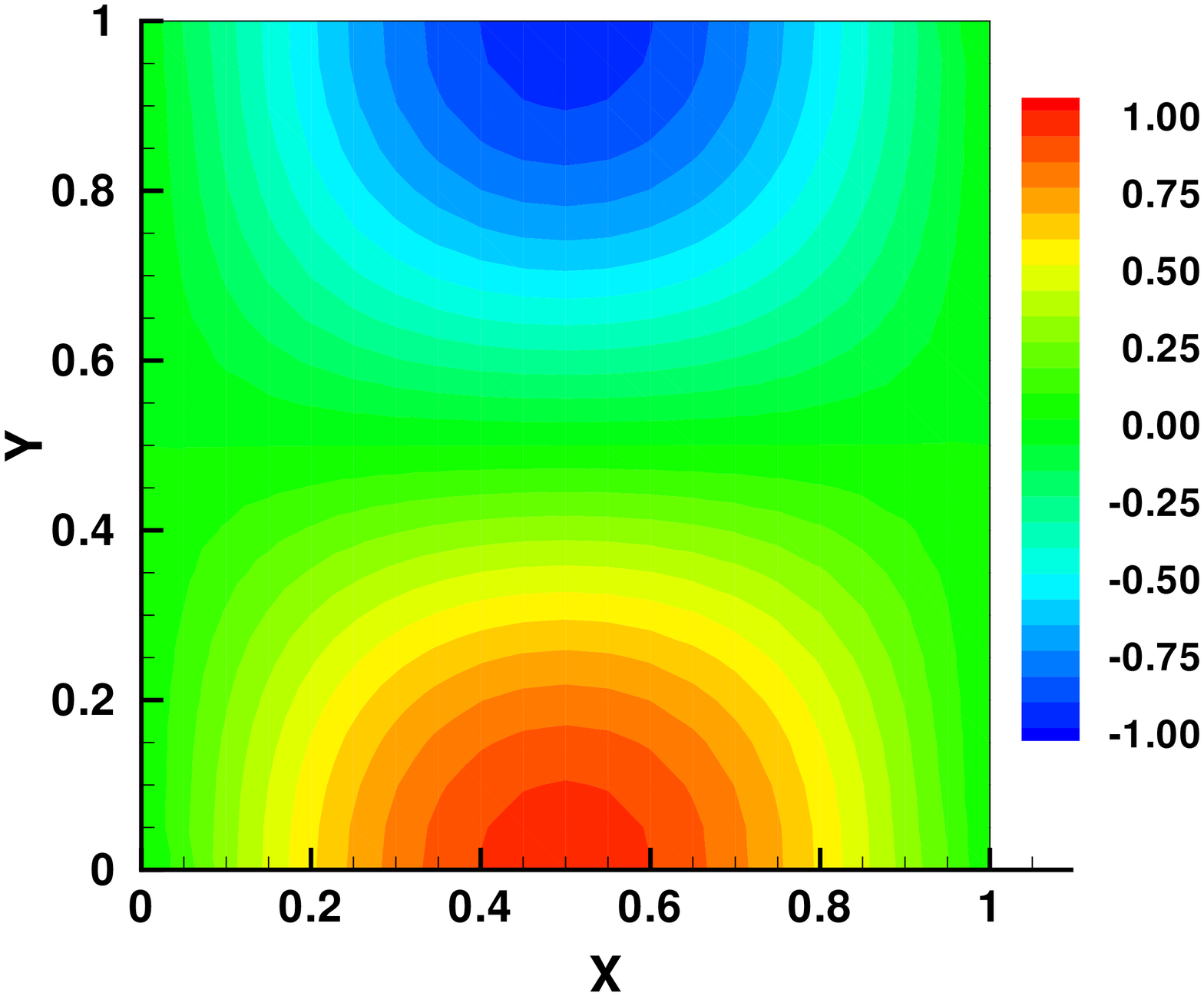}}
  \subfigure[$y$-velocity using Q4 elements]{
    \includegraphics[scale=0.3]{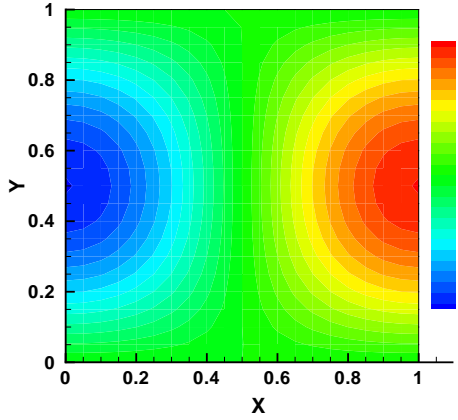}}
  \subfigure[$y$-velocity using T3 elements]{
    \includegraphics[scale=0.3]{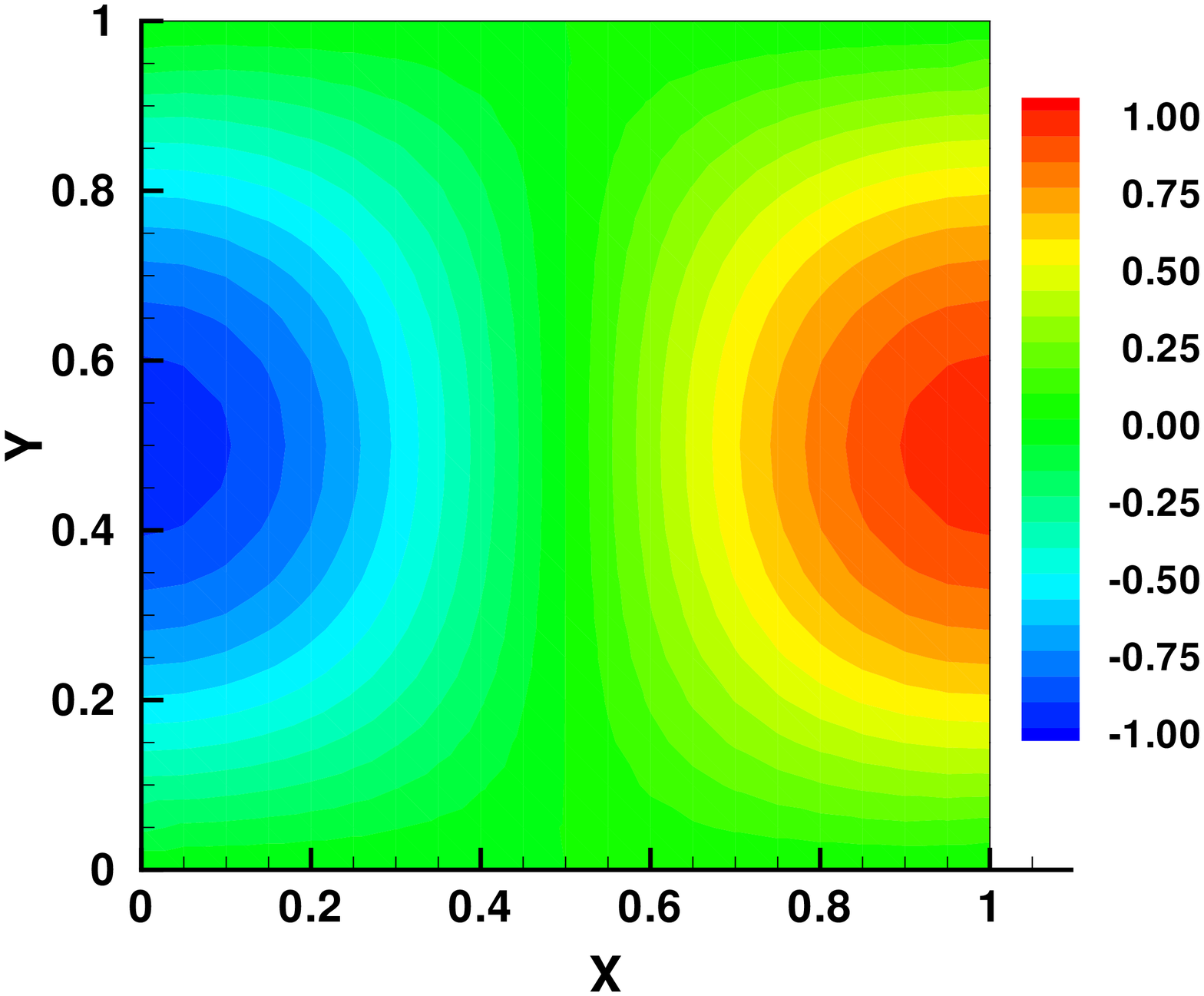}}
  \caption{Structured meshes: contours of $x$-velocity (top) and $y$-velocity (bottom) 
    using $100$ uniform four-node quadrilateral (Q4) elements and $200$ three-node 
    triangular (T3) elements. \label{Fig:MDarcy_2D_Problem_convergence_V}}
\end{figure}

\begin{figure}
  \centering
  \subfigure[Convergence for four-node quadrilateral elements]{
    \includegraphics[scale=0.3]{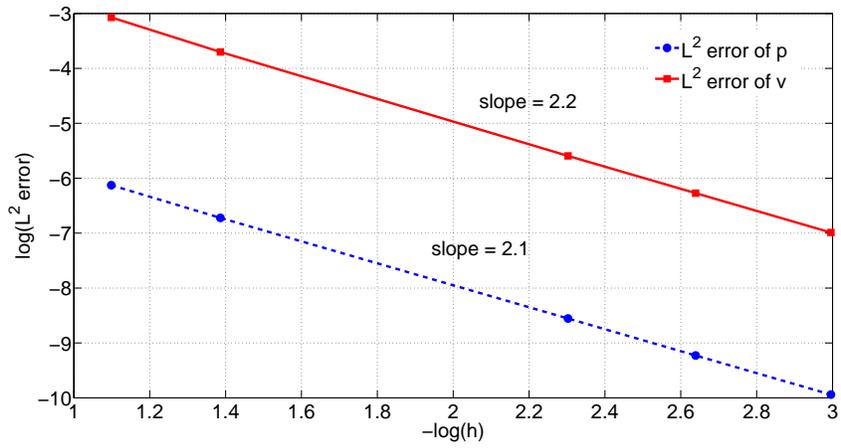}}
  \subfigure[Convergence for three-node triangular elements]{
    \includegraphics[scale=0.3]{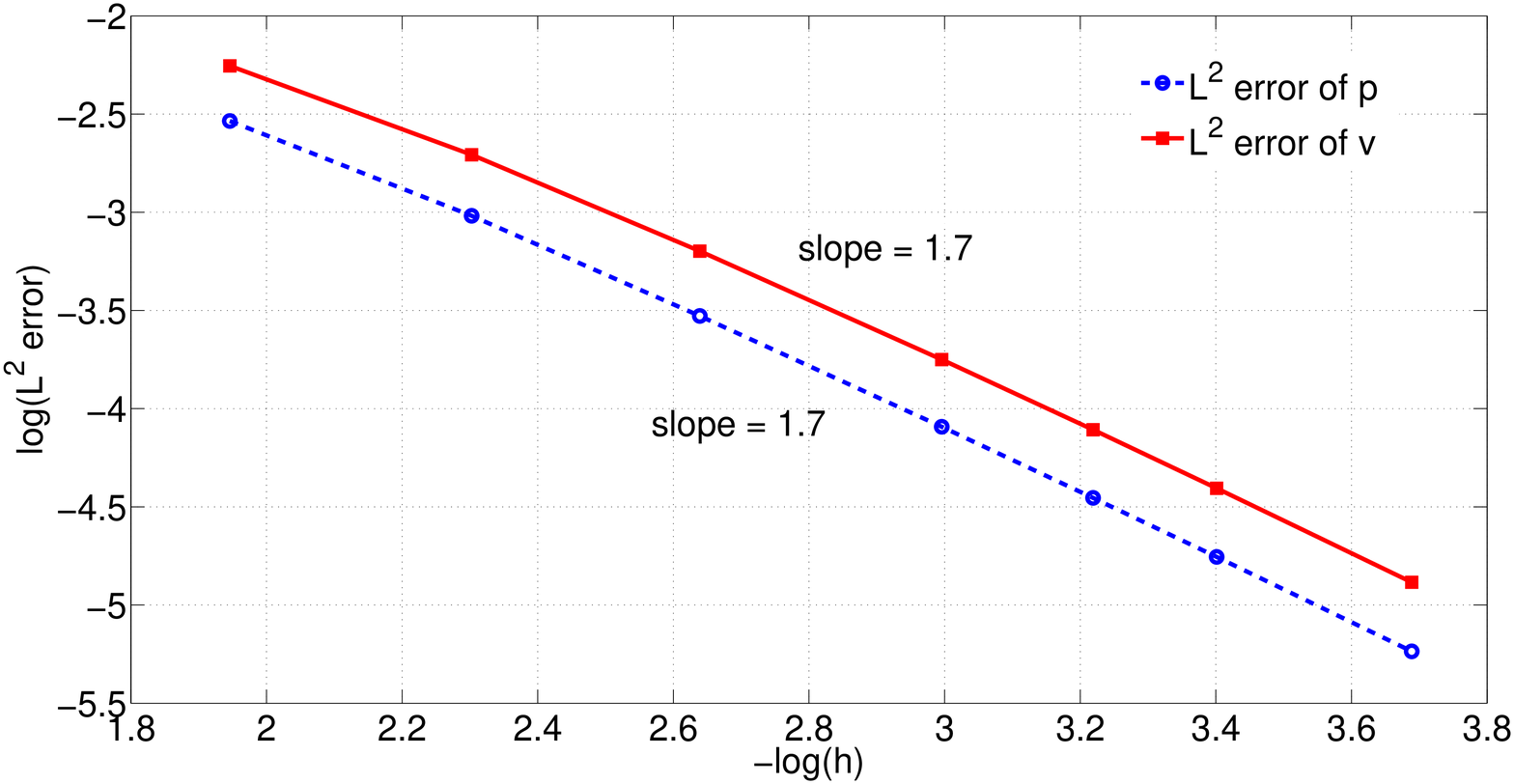}}
  \caption{2D $h$-convergence analysis: Rate of convergence of pressure and velocity in 
    the $L^{2}(\Omega)$ norm for four-node quadrilateral and three-node triangular 
    finite elements. \label{Fig:MDarcy_rate_of_h_convergence}}
\end{figure}

\begin{figure}
  \centering
  \subfigure[Computational mesh]{
    \includegraphics[scale=0.3]{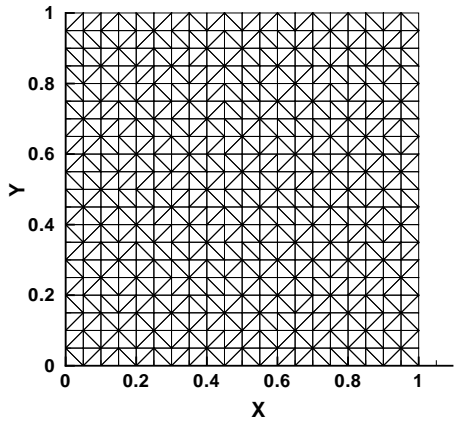}}
  \subfigure[Pressure contours]{
    \includegraphics[scale=0.3]{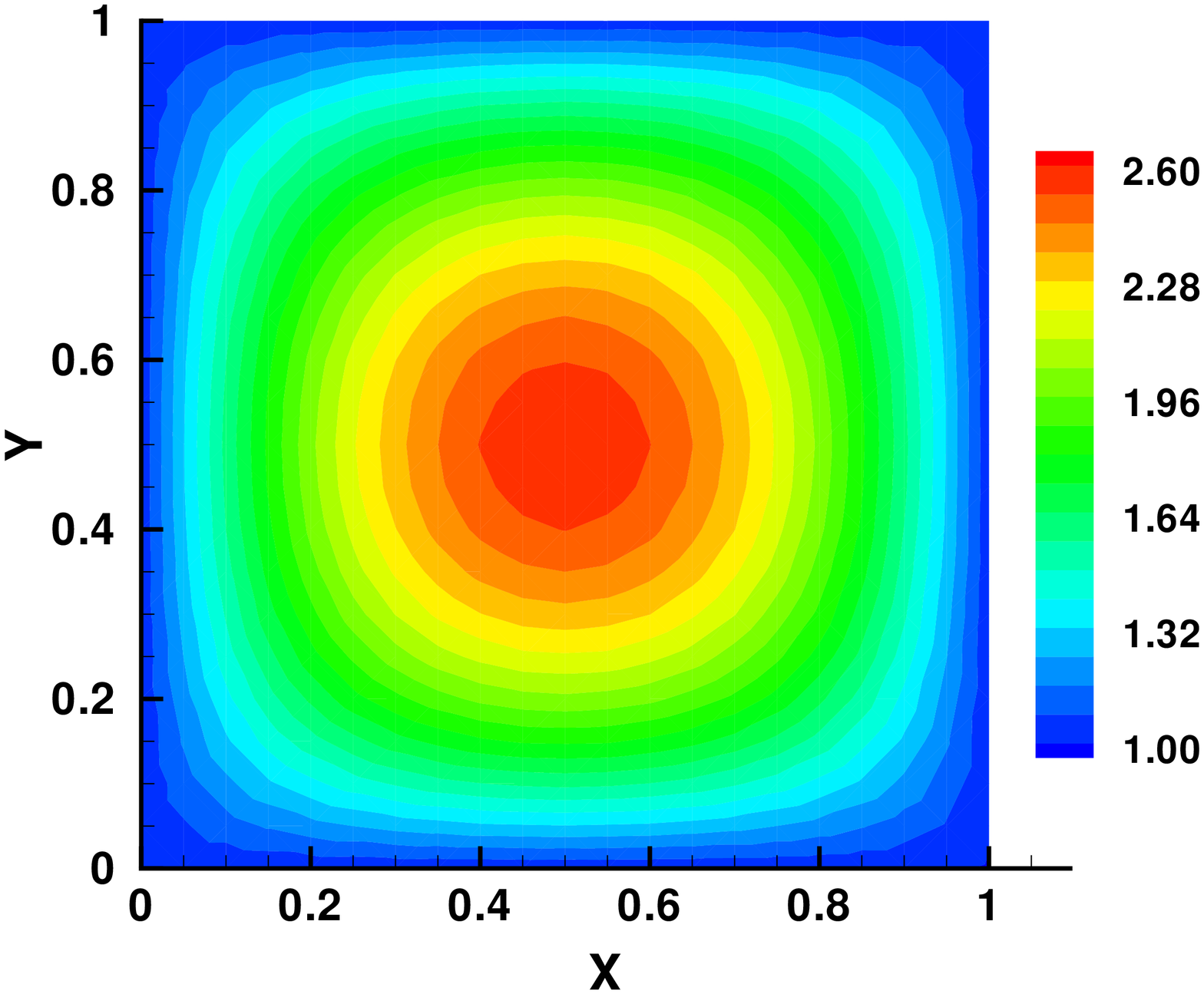}}
  \subfigure[$x$-velocity contours]{
    \includegraphics[scale=0.3]{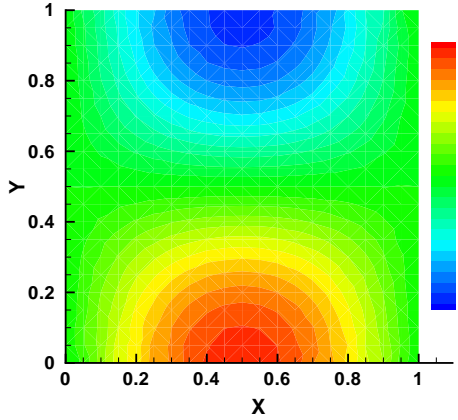}}
  \subfigure[$y$-velocity contours]{
    \includegraphics[scale=0.3]{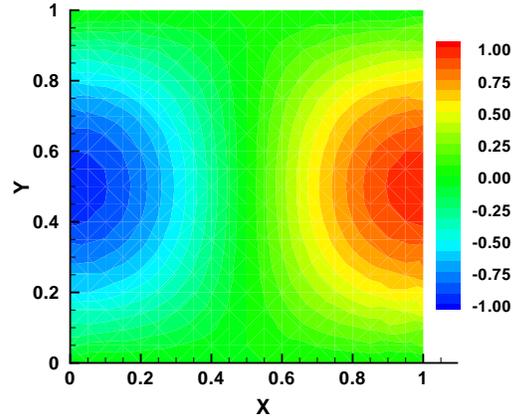}}
  \caption{Unstructured three-node triangular mesh: mesh (top-left), pressure (top-right), 
    $x$-velocity (bottom-left) and $y$-velocity (bottom-right). Barus' formula is employed 
    with $\alpha_0 = 1$ and $\beta = 2$. \label{Fig:MDarcy_2D_problem_Delaunay}}
\end{figure}

The contours of pressure and velocity are shown in Figures \ref{Fig:MDarcy_2D_Problem_convergence_Pre} 
and \ref{Fig:MDarcy_2D_Problem_convergence_V}. The rates of $h$-convergence of these finite elements 
for the aforementioned problem are shown in Figure \ref{Fig:MDarcy_rate_of_h_convergence}. (We have 
used the natural logarithm in these figures.) For the chosen problem, the proposed formulation converges 
in both pressure and velocity with respect to $h$-refinement. 
The above problem \eqref{Eqn:2D_problem_body_force} is also solved using an unstructured triangular mesh. 
The obtained numerical results and mesh are shown in Figure \ref{Fig:MDarcy_2D_problem_Delaunay}, quadratic 
convergence of the Newton-Raphson solution procedure for this problem is illustrated in Table 
\ref{Table:MDarcy_NR_convergence_Delaunay}, and the proposed mixed stabilized formulation performed well. 

\begin{table}
  \caption{Unstructured three-node triangular mesh: Convergence of Newton-Raphson solution 
    procedure for solving modified Darcy's equation. We have used Barus' formula with 
    $\alpha_0 = 1$ and $\beta = 2$. As one can see from the table, the solution procedure 
    exhibits terminal quadratic convergence. \label{Table:MDarcy_NR_convergence_Delaunay}}
  \begin{tabular}{cl} \\ \hline
    Iteration \# &  Norm of residual \\ \hline     
    0  &   1.353515824       \\
    1  &   3.594631185       \\
    2  &   1.781220654       \\ 
    3  &   31.25838881       \\
    4  &   0.169839581       \\
    5  &   0.015076303       \\
    6  &   3.125925098e-05   \\
    7  &   7.498072058e-10   \\
    8  &   6.532312307e-16   \\ \hline
  \end{tabular}
\end{table}

\subsection{A typical reservoir simulation}
\begin{figure}
  \subfigure{
    \includegraphics[scale=0.5]{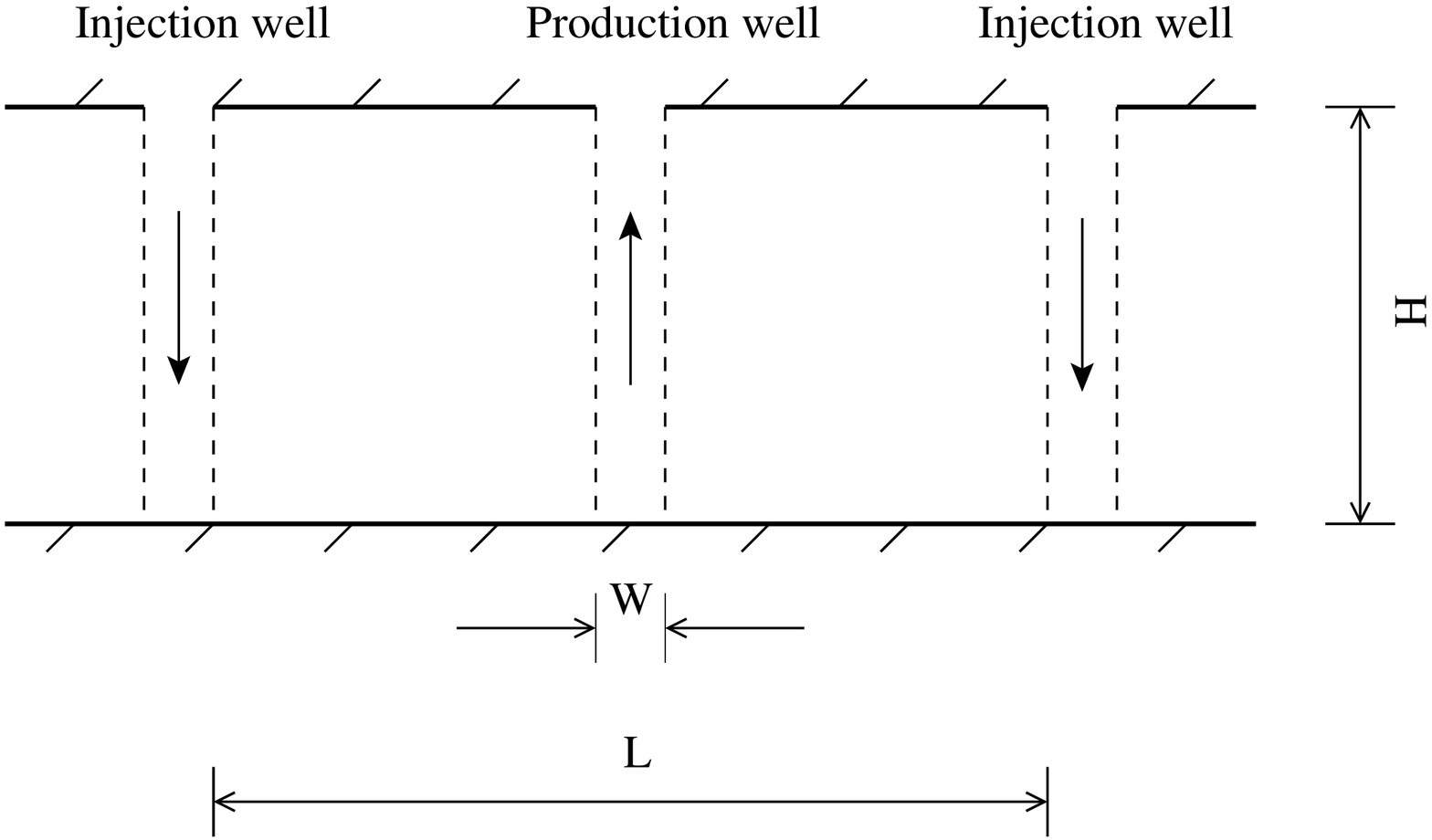}}
  \subfigure{
    \psfrag{vn}{$\boldsymbol{v}\cdot\boldsymbol{n} = 0$}
    \psfrag{g}{$\Gamma^{\mathrm{N}}$}
    \psfrag{D1}{$\Gamma^{\mathrm{D}}_1$}
    \psfrag{D2}{$\Gamma^{\mathrm{D}}_2$}
    \psfrag{p1}{$p(\boldsymbol{x}) = p_{\mathrm{atm}}$}
    \psfrag{p2}{$p(\boldsymbol{x}) = p_{\mathrm{enh}}$}
    \includegraphics[scale=0.5]{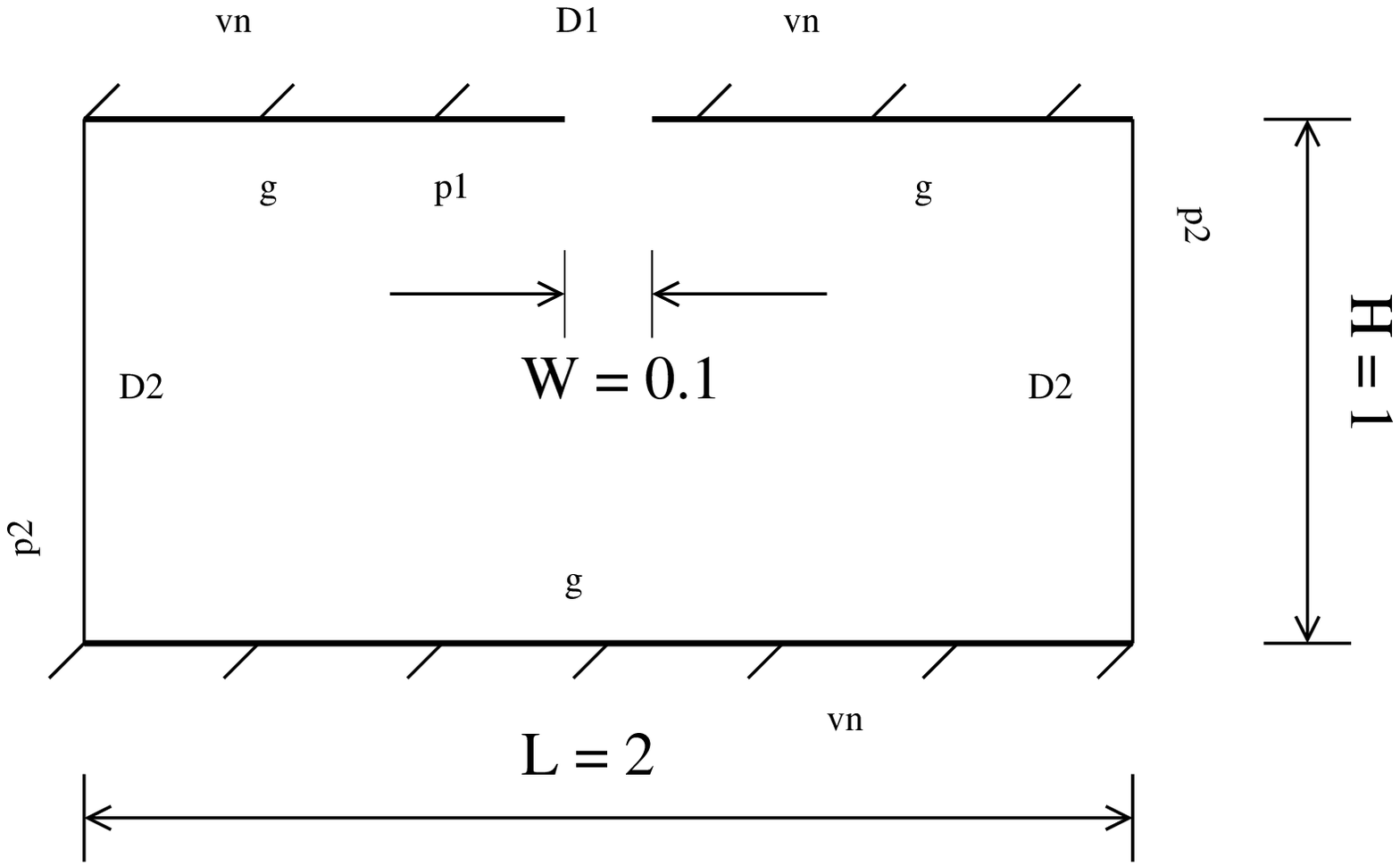}}
  \caption{Top figure is a pictorial description of a reservoir. Water, steam or Carbon dioxide 
    is pumped through injection wells, and crude oil is collected at the production well. Bottom 
    figure shows the computational idealization of the reservoir with appropriate boundary 
    conditions. On $\Gamma^{\mathrm{D}}_1$ we prescribe $p(\boldsymbol{x}) = p_{\mathrm{atm}}$, 
    and on $\Gamma^{\mathrm{D}}_2$ we prescribe $p(\boldsymbol{x}) = p_{\mathrm{enh}}$. 
    \label{Fig:MDarcy_reservoir}}
\end{figure}

\begin{figure}
  \includegraphics[scale=0.6]{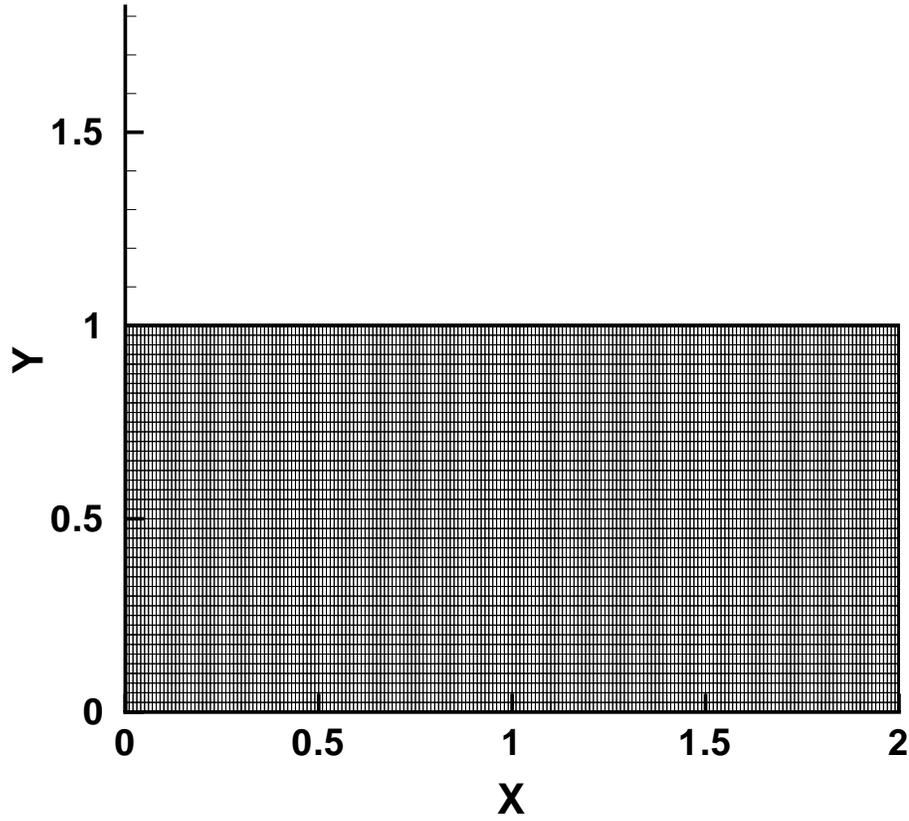}
  \caption{A four-node quadrilateral (structured) finite element mesh used in the reservoir 
    simulation. There are $8000$ elements, $8241$ nodes, and each node has three degrees-of-freedom 
    ($x$- and $y$-velocities, and pressure). \label{Fig:MDarcy_Q4_mesh}}
\end{figure}

Figure \ref{Fig:MDarcy_LLNL_EOR} shows how enhanced oil recovery is achieved in the field. 
A corresponding computational idealization of enhanced oil recovery is shown in Figure 
\ref{Fig:MDarcy_reservoir}, which will be discretized using finite elements. The computational 
mesh is shown in Figure \ref{Fig:MDarcy_Q4_mesh}. It should be noted that the corner points at 
the production well are reentrant corners, and hence the velocity $\boldsymbol{v}$ (which is 
proportional to gradient of pressure) is infinity at those reentrant corner points. We now 
compare the solutions obtained using Darcy's equation (constant drag function) and modified 
Darcy's equation using Barus' formula. 

For this problem, we have taken $\bar{\alpha}_0 = 1$, $\bar{\beta} = 0.005$, $P = 1 \; \mathrm{atm}$, 
$B = 9.81 \; \mathrm{m/s^2}$, $\bar{\boldsymbol{b}} = [0, -1]^{T}$, $\mathcal{A} = 1$ and $\mathcal{C} = 1$. 
The top surface at the production well is at atmospheric pressure (see Figure \ref{Fig:MDarcy_reservoir}), 
which is given by the non-dimensional parameter $\bar{p}_{\mathrm{atm}} = 1$ (as we have taken the reference 
pressure to be $P = 1 \; \mathrm{atm}$).

\begin{figure}
  \subfigure{
    \includegraphics[scale=0.5]{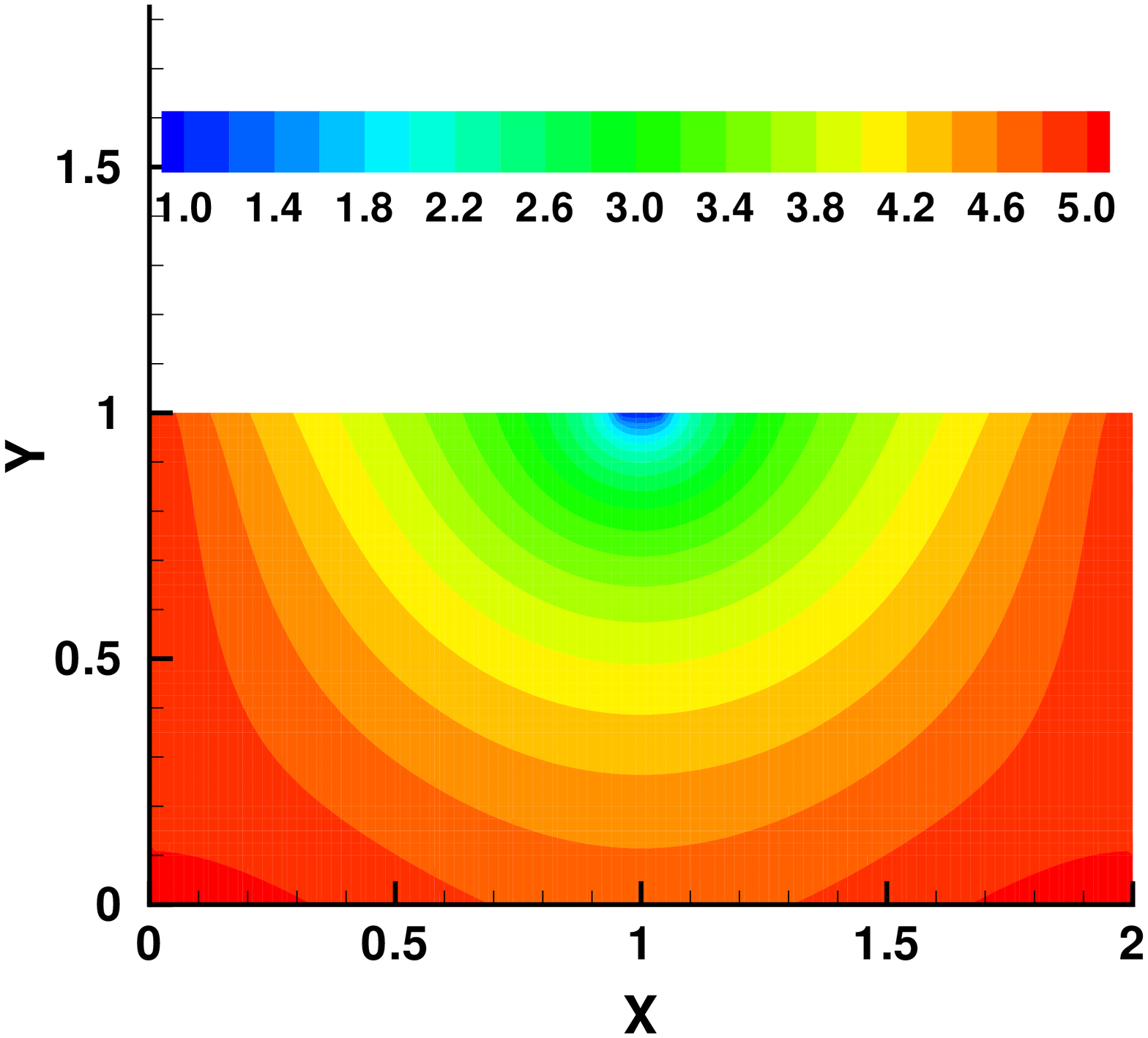}}
  \subfigure{
    \includegraphics[scale=0.5]{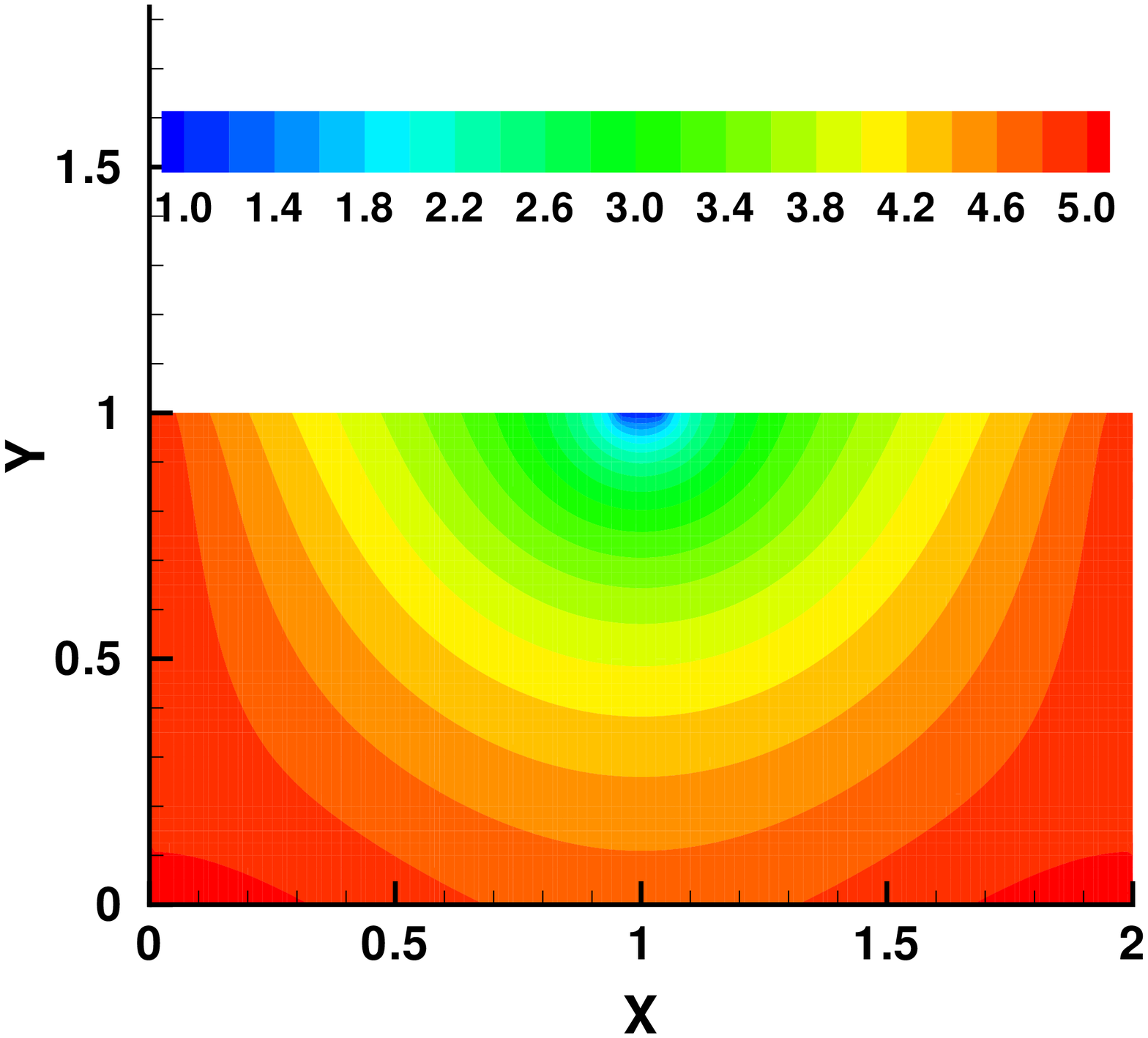}}
  \caption{Pressure contours using Darcy's equation (top) and modified Darcy's equation 
    (bottom) for $\bar{p}_{\mathrm{enh}} = 5$. \label{Fig:MDarcy_well_p_5}}
\end{figure}

\begin{figure}
  \subfigure{
    \includegraphics[scale=0.5]{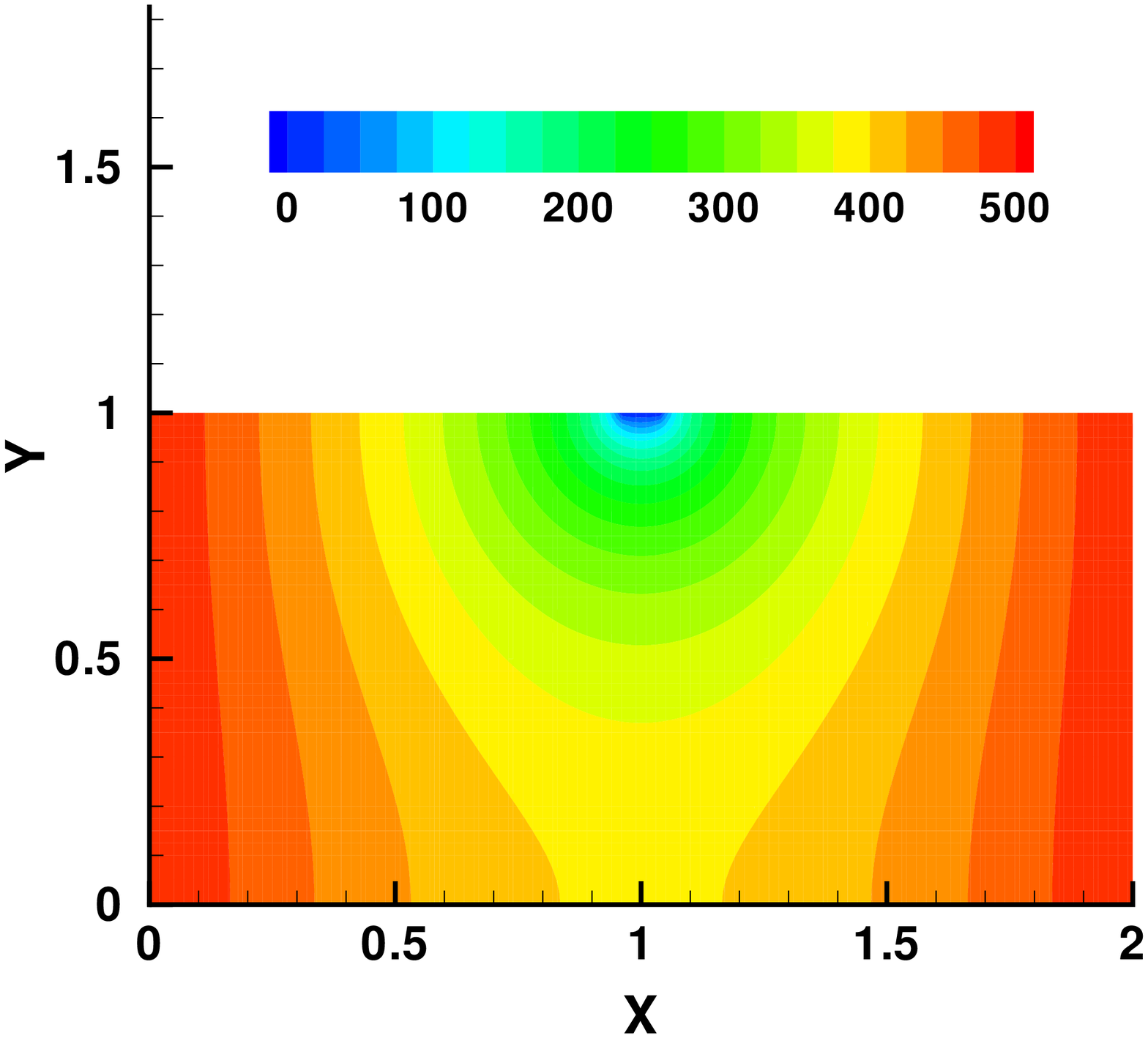}}
  \subfigure{
    \includegraphics[scale=0.5]{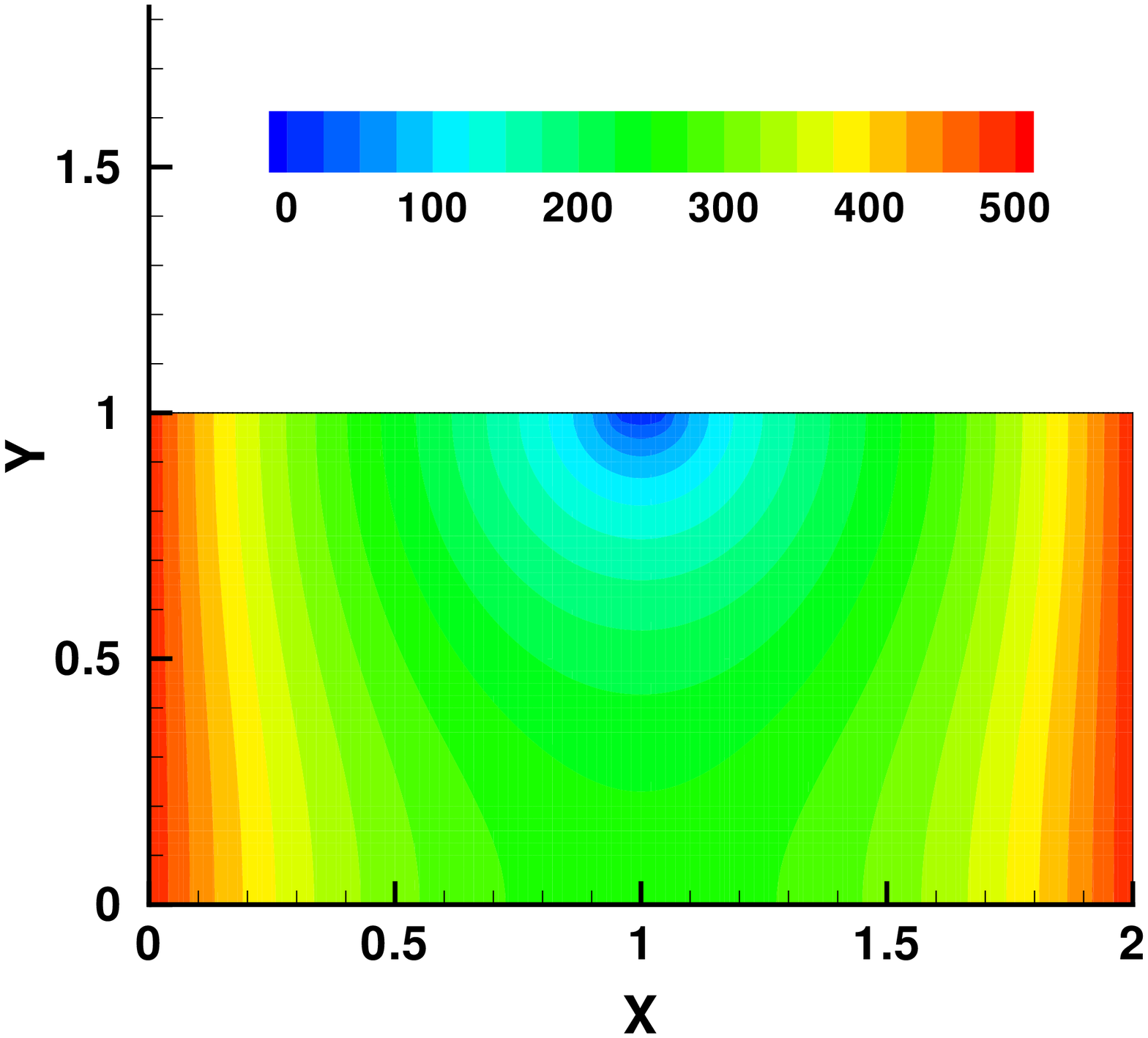}}
  \caption{Pressure contours using Darcy's equation (top) and modified Darcy's equation 
    (bottom) for $\bar{p}_{\mathrm{enh}} = 500$. \label{Fig:MDarcy_well_p_500}}
\end{figure}

In Figures \ref{Fig:MDarcy_well_p_5} and \ref{Fig:MDarcy_well_p_500}, pressure contours are shown for the 
cases $\bar{p}_{\mathrm{enh}} = 5$ and $\bar{p}_{\mathrm{enh}} = 500$, respectively. For low prescribed 
pressures at the injection wells, the pressure profiles obtained using Darcy's equation and modified 
Darcy's equation are similar. However, this is not the case for higher pressures at the injection wells, 
and the pressure profiles using Darcy's and modified Darcy's equations are significantly different. The 
pressure using modified Darcy's equation exhibit steeper gradients near the injection wells. It should 
be noted that pressure gradients play an important role in the study of fracture of porous medium. 

\begin{figure}
  \centering
  \includegraphics[scale=0.4]{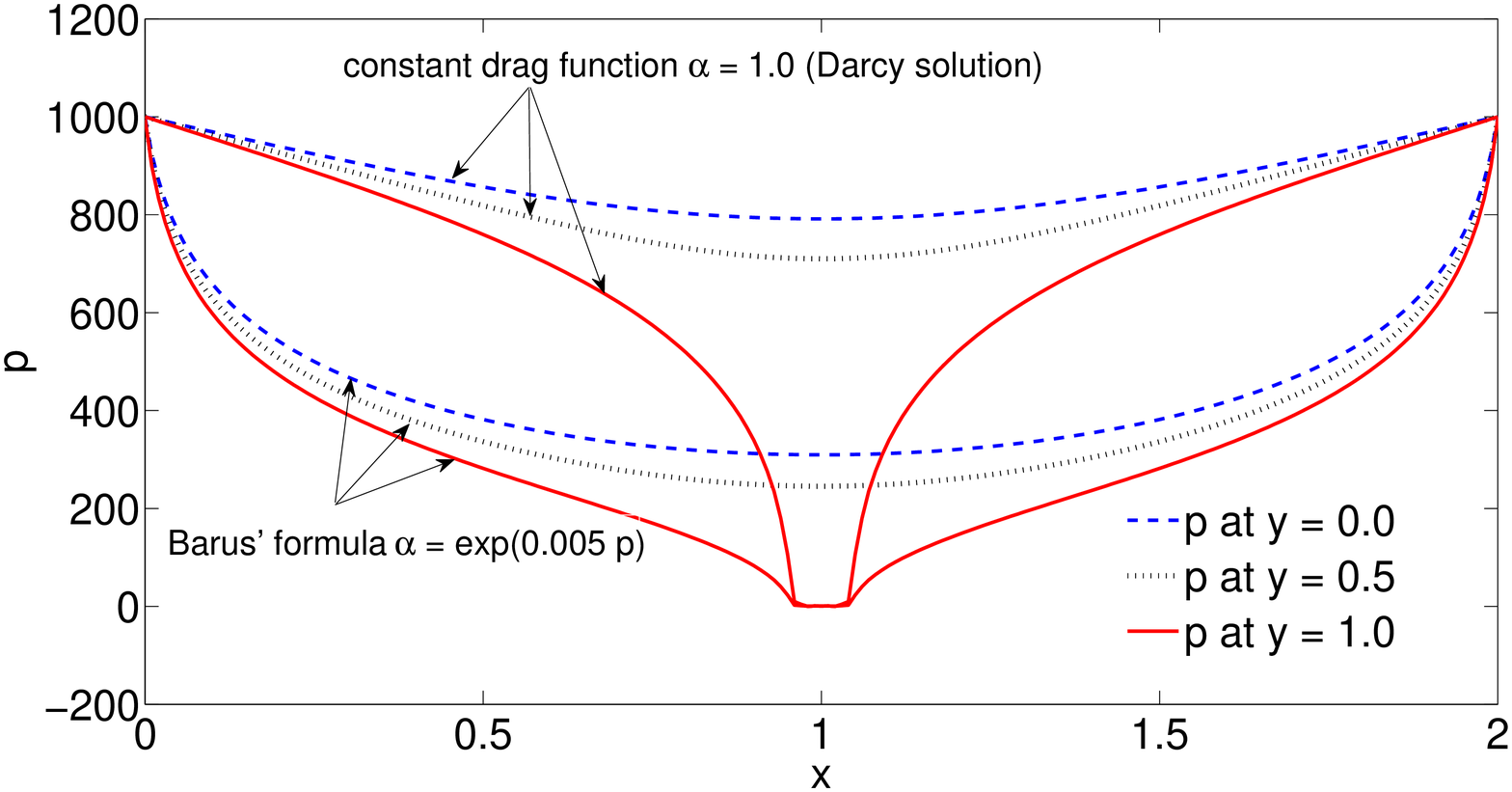}
  \caption{Variation of pressure with respect to horizontal distance from the production 
    well at various depths of the reservoir using Darcy's equation and modified Darcy's 
    equation for $\bar{p}_{\mathrm{enh}} = 1000$. \label{Fig:MDarcy_pressure_variation}}
\end{figure}

In Figure \ref{Fig:MDarcy_pressure_variation}, the variations of pressure with respect to the horizontal 
distance from the production well at various depths are plotted for Darcy's and modified Darcy's equations. 
In this study we have taken $\bar{p}_{\mathrm{enh}} = 1000$. As one can see from Figure 
\ref{Fig:MDarcy_pressure_variation}, the pressures from Darcy's equation and modified Darcy's equation are 
qualitatively and quantitatively different. Most importantly, the cone shapes are completely different. 
(Cone shape is the variation of pressure with respect to the horizontal distance from the production well, 
and is considered an important parameter in assessing the well performance and also regulating reservoir-operating 
conditions.) For the case of Darcy's equation, the graph is concave upwards, while the graph using modified Darcy's 
equation (with Barus' formula) is convex upwards. 

\begin{figure}
  \centering
  \includegraphics[scale=0.4]{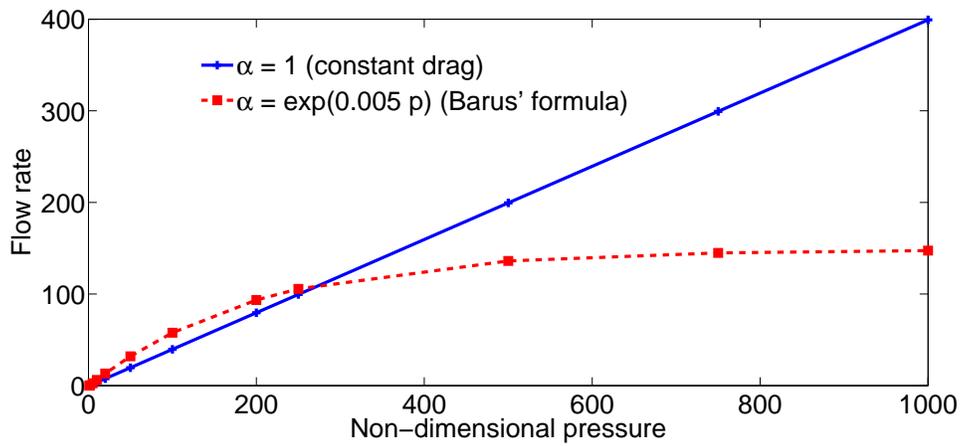}
  \caption{Total flow rate at the production well using Darcy's equation and modified Darcy's 
    equation (using Barus' formula) for various values of $\bar{p}_{\mathrm{enh}}$. 
    \label{Fig:MDarcy_total_flux}}
\end{figure}

In Figure \ref{Fig:MDarcy_total_flux} we compared the total flux at the production well for 
various pressures at the injection wells (which is given by $\bar{p}_{\mathrm{enh}}$) using 
Darcy's equation and modified Darcy's equation (based on Barus' formula). Even for this case 
the pressure at the production well is $\bar{p}_{\mathrm{atm}} = 1$. The total flux at the 
production well is calculated using 
\begin{align}
  \int_{\Gamma^{\mathrm{D}}_1} \boldsymbol{v} \cdot \boldsymbol{n} \; \mathrm{d} \Gamma
\end{align}
where the part of boundary $\Gamma^{\mathrm{D}}_1$ is indicated in Figure \ref{Fig:MDarcy_reservoir}. 
As one can see in Figure \ref{Fig:MDarcy_total_flux}, modified Darcy's equation predicts a ceiling 
for the total flux with respect to pressure, which is what is expected physically. On the other hand, 
Darcy's equation predicts continued linear increase of the total flux with respect to the pressure. 

\emph{
One can expect even more interesting departures from the classical results when the Darcy-Forchheimer 
and Darcy-Forchheimer-Brinkman equations are modified to incorporate the effect of pressure on the 
viscosity and drag function. Studying the interactions between the solid deformation and fluid flow, 
with the effects due to the pressure being taken into account will provide even more insights into 
the real problem, and is part of our future work.}

%
\section{CONCLUSIONS}
\label{Sec:MDarcy_Conclusions}
In this paper we have studied the flow of an incompressible fluid through a rigid 
porous solid under high pressure and pressure gradients. For such problems, Darcy's 
equation is not a good model as it assumes constant viscosity and drag coefficient, 
which is contrary to experimental evidence. Here, we have considered a modification 
to Darcy's equation wherein the drag coefficient is a function of pressure, which is 
the case with many liquids for large range of pressures. We have allowed the drag 
coefficient to vary with pressure in different ways, which are primarily motivated 
by experimental observations. 
We have also presented a new stabilized mixed formulation for the modified Darcy's 
equation along with the incompressibility constraint. The proposed formulation allows 
equal-order interpolation for pressure and velocity, which is not stable under the 
classical mixed formulation. A noteworthy feature of the formulation is that the 
stabilization parameter neither involves mesh-dependent nor material parameters. 
We have presented representative numerical results to show that the proposed formulation 
performs well. We have also shown that the solution (both velocity and pressure fields) 
based on the modified Darcy's equation is significantly different from the solution 
based on Darcy's equation. In particular, using a representative reservoir simulation, 
we have shown that the Darcy's model may predict significantly higher flow rate at high 
pressure as it does not account for variation of viscosity with respect to pressure.

%
\section{APPENDIX}
\label{Sec:MDarcy_Appendix}
\subsection{Notation and definitions}
We now define some quantities that will be useful in writing the finite element residual vector and 
tangent stiffness matrix in a compact manner. Let $\boldsymbol{A}$ and $\boldsymbol{B}$ be matrices 
of size $n \times m$ and $p \times q$, respectively. That is,
\begin{align*}
\boldsymbol{A} = \left[ \begin{array}{ccc} a_{1,1} & \hdots & a_{1,m}  \\
\vdots & \ddots & \vdots \\
a_{n,1} & \hdots & a_{n,m} \end{array} \right] \ ; \ \boldsymbol{B} = 
\left[ \begin{array}{ccc} b_{1,1} & \hdots & b_{1,q}  \\
\vdots & \ddots & \vdots \\
b_{p,1} & \hdots & b_{p,q} \end{array} \right]
\end{align*}
The \emph{Kronecker product} of these matrices is an $np \times mq$ matrix, and is defined as
\begin{align*}
\boldsymbol{A} \odot \boldsymbol{B} := \left[ \begin{array}{ccc} a_{1,1}\boldsymbol{B} & 
\hdots & a_{1,m}\boldsymbol{B}  \\ \vdots & \ddots & \vdots \\
a_{n,1}\boldsymbol{B} & \hdots & a_{n,m}\boldsymbol{B} \end{array} \right]
\end{align*}
Note that the Kronecker product can defined for \emph{any} two given matrices (irrespective of their 
dimensions). The $\mathrm{vec}[\cdot]$ operator is defined as
\begin{align*}
\mathrm{vec}[\boldsymbol{A}] := \left[ \begin{array}{c} a_{1,1} \\
\vdots \\
a_{1,m} \\
\vdots \\
a_{n,1} \\
\vdots \\
a_{n,m} \end{array} \right]
\end{align*}
Some relevant properties of Kronecker product and $\mathrm{vec}[\cdot]$ operator are as follows: 
\begin{itemize}
\item $\mathrm{vec}[\boldsymbol{A C B}] = \left(\boldsymbol{B}^{T} \odot \boldsymbol{A}\right) 
  \mathrm{vec}[\boldsymbol{C}]$
\item $\left(\boldsymbol{A} \odot \boldsymbol{B} \right) \left(\boldsymbol{C} \odot \boldsymbol{D} \right) = 
  \left(\boldsymbol{A} \boldsymbol{C} \odot \boldsymbol{B} \boldsymbol{D} \right)$
\item $\mathrm{vec}[\boldsymbol{A} + \boldsymbol{B}] = \mathrm{vec}[\boldsymbol{A}] + 
  \mathrm{vec}[\boldsymbol{B}]$
\end{itemize}
For a detailed discussion on Kronecker products and $\mathrm{vec}[\cdot]$ operator see References 
\cite{Graham_Kronecker,Laub,Nicholson}.

\bibliographystyle{unsrt}
\bibliography{../Master_References/Master_References,../Master_References/Books}
\newpage

\clearpage

\end{document}